\documentclass[a4paper,12pt,final]{amsart}
\usepackage[utf8]{inputenc}
\usepackage{eurosym}
\usepackage{amsfonts}
\usepackage{amsmath}
\usepackage{graphicx}
\usepackage{graphics}
\usepackage{epsfig}
\usepackage{epstopdf}
\usepackage{geometry}
\usepackage{graphics}
\usepackage{color}
\usepackage{caption}
\usepackage{ulem}
\usepackage{showkeys}

\begin{document}

\title[]{Numerical study of the Serre-Green-Naghdi equations in 2D}

\author{Sergey Gavrilyuk}
\address{Aix-Marseille Universit\'e and CNRS UMR 7343 IUSTI, 5 rue Enrico Fermi, 13453 Marseille, France\\
E-mail sergey.gavrilyuk@univ-amu.fr}
\author{Christian Klein}
\address{Institut de Math\'ematiques de Bourgogne, UMR 5584\\
                Universit\'e de Bourgogne, 9 avenue Alain Savary, 21078 Dijon
                Cedex, France\\
				Institut Universitaire de France\\
    E-mail Christian.Klein@u-bourgogne.fr}

\date{\today}
\begin{abstract}

A detailed numerical  study of solutions to the Serre-Green-Naghdi (SGN) 
equations in 2D with vanishing curl of the velocity field is presented. 
The transverse stability of line 
solitary waves, 1D solitary waves being exact solutions of the 2D 
equations independent of the second variable, is established numerically. 
The study of localized 
initial data as well as crossing 1D solitary waves does not give an 
indication  of existence  of stable structures in SGN solutions localized in two 
spatial dimensions. For the numerical experiments, an approach based on a 
Fourier spectral method with a Krylov subspace technique is 
applied. 
\end{abstract}

\thanks{CK thanks for  support by   the EIPHI Graduate School (contract ANR-17-EURE-0002) and by the 
European Union Horizon 2020 research and innovation program under the 
Marie Sklodowska-Curie RISE 2017 grant agreement no. 778010 IPaDEGAN. }

\maketitle

\section{Introduction}
We present a numerical study of solutions to the 
Serre--Green--Naghdi (SGN) equations. The SGN equations are classical dispersive shallow water equations which can be derived from  the free surface Euler equations both   by depth averaging  \cite{Serre_53,Su_Gardner_1969,Green_74,Green_76} 
	and via Hamilton's  principle of stationary action~\cite{Salmon,Salmon_1998}. 
	In dimensional form, the SGN equations for a flat bottom are: 
	\begin{subequations}
		\label{eq:SGN}
		\begin{align}
			& h_{t}+\nabla\cdot(h\overline{\mathbf u})=0,\label{eq:SGN-mass} \\
			& (h\overline{\mathbf u})_{t} +
		\nabla\cdot\left(h \overline{\mathbf 
			u}\otimes\overline{\mathbf u}+p{\mathbf I}\right)= 0, \quad p=\frac{gh^2}{2}+
		\frac{ h^2}{3} \ddot h . \label{eq:SGN-mom}
			\end{align}
			\end{subequations}
				Here $h$ is the water depth, $\overline{\mathbf 
				u}=(u_x,u_y)^T$ is the depth averaged velocity, $\mathbf x=(x,y)^T$ are the classical Cartesian coordinates, $g$ is the gravity acceleration,
	$\nabla\cdot$ is the  divergence operator,  $\displaystyle \dot 
		h=\frac{\partial h}{\partial t}+\overline{\mathbf u}\cdot \nabla h$ is the material derivative along the depth averaged velocity, 
		two ``dots'' denote the corresponding second material derivative. 
		As a consequence, they admit the energy conservation law in 
		the form: 
		\begin{equation}
			\label{energy}
				\left(hE\right)_{t}+
			\nabla\cdot\left(hE\overline{\mathbf u}+
		p\overline{\mathbf u}\right)
			= 0, \quad E=\frac{\vert\overline{\mathbf u}\vert^2}{2}+\frac{gh}{2}+\frac{\dot{h}^2}{6}.
			\end{equation}
			
A class of generalized potential flows which are 
exact solutions to the 2D SGN equations was  introduced in  \cite{Gavrilyuk_Teshukov_2001}. In 
particular, this class of flows is well defined if  the 
non-cavitation condition ($h>0$) is fulfilled. It is assumed 
throughout the manuscript.  To study such 2D solutions 
numerically,        
we apply an extension of the numerical approach \cite{DK} for the 1D 
SGN equation to two spatial dimensions. The method is again based on 
a Fourier spectral method in the spatial variables with the Krylov 
subspace technique GMRES \cite{GMRES} to invert an elliptic equation 
for an auxiliary quantity (the appearence of the quantity $\ddot{h}$ 
in (\ref{eq:SGN-mom}) implies that at least one elliptic equation has 
to be solved in each time step as will be detailed below). The time 
integration will be done with the standard explicit 4th order 
Runge-Kutta method. 

This code will be applied to interesting initial data in the context 
of the SGN equations. First we address the transverse stability of 
the line solitary waves: the explicitly known solitary wave solution 
for the 1D SGN equations (\ref{sol}) can be seen as a $y$-independent 
exact solution to the 2D SGN equations, i.e., an infinitely extended 
travelling wave called \textit{line solitary wave}. We study various
perturbations of this line solitary wave numerically for a range of 
velocities relevant in applications. In this sense we study what is 
known as the transverse stability of such infinitely extended 
structures. The numerical results give strong evidence to the 
following\\
\textbf{Main conjecture I:}\\
\textit{The line solitary wave (\ref{sol}) is transversely stable as a 
solution to the 2D SGN equations. }

The stability of the line solitary wave is an indication that the 
behavior of solutions to the SGN equations is similar to the one of 
solutions to the Kadomtsev-Petviashvili (KP) II equations, see 
\cite{KSbook} for a recent comprehensive review with many references. 
For the KP II equation, the Korteweg-de Vries (KdV) soliton gives a 
line soliton solution which is stable for KP II, but not for the 
so-called KP I equation. In the latter case the line soliton is 
unstable against the formation of so-called lump solitons, travelling 
waves localised in two dimensions. The KP II equation does not have 
travelling wave solutions localised in two dimensions, and it has a 
\textit{defocusing} effect on localised initial data: there are no 
stable structures in their solutions localised in 2D, hump-like 
initial data are simply radiated away to infinity. By studying 
localised initial data for SGN, we give strong evidence for the fact 
that SGN has a defocusing effect as KP II, and that there are no 
stable solitary waves localised in two dimensions. We get the 
following\\
\textbf{Main conjecture II:}\\
\textit{There are no stable SGN solitary waves localised in two 
dimensions. The equation has a defocusing effect.  }

The paper is organised as follows: in section 2 we collect some basic 
facts of the 2D SGN equations. In section 3 we briefly describe the 
applied numerical approach. The transverse stability of the line 
solitary waves is studied in section 4. In section 5 we consider the 
time evolution of the formal superposition of two line solitary 
waves, one in $x$, the other in $y$. Again no stable structures are 
observed in the time evolution of these data. In section 6 we study 
the time evolution of hump-like initial data. It is shown that 
Gaussian initial data evolve into an annular structure with a central 
depression. We add some concluding remarks in section 7.

\section{Basic facts}
In this section, we collect some basic facts on the SGN equation and 
give a form of the equations suitable for the planned numerical 
treatment.

\subsection{Introduction of the potential $\varphi$}

One can  introduce the  variable  ${\mathbf v}=(v_x,v_y)^T$  defined as 
 \cite{Salmon,Salmon_1998,Gavrilyuk_Teshukov_2001}:
 \begin{equation}
{\mathbf v}=\overline{\mathbf u}-\frac{1}{3h}{\nabla}(h^3\;\nabla\cdot\overline{\mathbf u})=\overline{\mathbf u}+\frac{1}{h}{\nabla}\left(\frac{h^2\dot h}{3}\right)
	\label{K_1}
\end{equation}
and the potential $\varphi $  such that $\mathbf{v}=\nabla \varphi$. 
The variable $\mathbf v$ is  the tangent velocity of the fluid at the 
free surface \cite{Gavrilyuk_2014}, and $\varphi$ is the generalized 
flow potential. One can prove that if, initially, ${\rm curl} \; 
\mathbf v=0$, it will stay zero for any time \cite{Gavrilyuk_Teshukov_2001}.  In the  following, we consider this special class of generalized potential flows.  The notion of generalized potential flows has a sense only in the 2D case because all one-dimensional flows are obviously irrotational. 
In the one-dimensional case the SGN equations admit, in particular,  
one-dimensional solitary wave solutions \cite{Su_Gardner_1969}. The 
linear  stability of solitary waves of small amplitude  has been 
analytically proven in \cite{Li_2001}, see for instance \cite{DK} for 
a numerical investigation. Here, we study their stability  with respect to multi-dimensional perturbations in the  class of generalized potential flows introduced above.

\subsection{Governing equations of  potential flows of the SGN equations}
Introducing the variable $\sigma$ 
\begin{equation}
	\sigma = \frac{h^2}{3}\dot{h}, 
\end{equation} 
and as above the  potential $\varphi$,
\begin{equation}
\nabla\varphi= \overline{\mathbf u}+ \frac{1}{h}\nabla\sigma, 
\end{equation} 
one can rewrite the  momentum equation \eqref{eq:SGN-mom} in the  form of the generalized Bernoulli equation \cite{Gavrilyuk_Teshukov_2001}:
	\begin{equation}
		\displaystyle \varphi_t+\nabla \varphi\cdot\overline{\mathbf 
		u}+gh-\frac{1}{2}\dot h^2-\frac{1}{2}\vert \overline{\mathbf u}\vert^2=gh_\infty, \quad \overline{\mathbf u}=\nabla \varphi-\frac{1}{h}\nabla \sigma.
\end{equation}
Here $h_{\infty}=const$ is the fluid depth at infinity, the  corresponding averaged velocity is vanishing at infinity. 
An equivalent form is  
	\begin{equation}
	\varphi_t+\frac{\vert \nabla\varphi\vert^2}{2}-\frac{\vert \nabla\sigma\vert^2}{2h^2}+gh-\frac{9}{2}\frac{\sigma^2}{h^4}=gh_\infty. 
\end{equation}
Complemented by the mass conservation law, the system reads:
\begin{subequations}
	\label{potential_SGN}
	\begin{align}
		&  h_t+\nabla\cdot(h\nabla\varphi)= \Delta \sigma,\label{1} \\
		& \varphi_t+\frac{\vert \nabla\varphi\vert^2}{2}-\frac{\vert \nabla\sigma\vert^2}{2h^2}+gh-\frac{9}{2}\frac{\sigma^2}{h^4}=gh_\infty, \label{2}\\
		&\frac{3\sigma}{h^3}-\nabla\cdot\left(\frac{\nabla \sigma}{h}\right)=-\Delta \varphi.\label{3}
	\end{align}
\end{subequations}
The equation \eqref{3} can be formally rewritten  in the form 
\begin{equation}
\mathcal{L}[\sigma]=-\Delta \varphi, \quad {\rm with} \quad \mathcal{L}[\sigma]=\frac{3\sigma}{h^3}-\nabla\cdot\left(\frac{\nabla \sigma}{h}\right).
\end{equation}
Then 
\begin{equation}
\sigma =-	\mathcal{L}^{-1}[\Delta\varphi].
\end{equation}
The boundary conditions as $\vert\mathbf x\vert\rightarrow\infty$ are:
\begin{equation}
\nabla\varphi\rightarrow 0, \quad \sigma \rightarrow 0, \quad h\rightarrow h_{\infty}>0. 
\label{bc}
	\end{equation}
These boundary conditions are adapted for the study of the solitary waves.

 \subsection{Conserved quantities}
The total mass, total momentum  and  total energy are conserved. Their  integral form for  boundary conditions \eqref{bc} is :
\begin{subequations}
	\begin{align}
		& \frac{d\mathcal{M}}{dt}=\frac{d}{dt}\iint_{-\infty}^{+\infty}\left(h-h_{\infty}\right)dxdy=0,
		\label{total_mass}	 \\
	& \frac{d\mathcal{P}}{dt}=\frac{d}{dt}\iint_{-\infty}^{+\infty}h \overline{\mathbf u}\;dxdy=0,
	\label{total_momentum}	\\
		& \frac{d\mathcal{E}}{dt}=\frac{d}{dt}\iint_{-\infty}^{+\infty}\left(\frac{h\vert\overline{\mathbf u}\vert^2}{2}+\frac{3\sigma^2}{2h^3}+\frac{g(h-h_{\infty})^2}{2}\right)dxdy=0.
		\label{total_energy}	
	\end{align}
\end{subequations}
Compared to \eqref{energy}, we used the total mass conservation  \eqref{total_mass} to transform the local energy density in \eqref{energy} to the form \eqref{total_energy}. 

\subsection{1D solitary waves}
Let $x$ be the one-dimensional space coordinate, and $u$ be the velocity in $x$--direction.  Let $\xi=x-ct$ be the travelling coordinate, with   $c>0$ being the wave velocity, and  $\overline{\mathbf{u}}=(u_c,0)^T$.  In the following, ``prime''  means the derivative with respect to $\xi$.
Then  the solution (for $c^2>gh_{\infty}$) is \cite{Su_Gardner_1969}:
\begin{equation} \label{sol}
h_{c}=h_{ \infty}+\frac{(c^2-gh_{ \infty})}{g}\mbox{sech}^2\left( 
\frac{\xi}{2ch_{ \infty}}\sqrt{3(c^2-gh_{ \infty})}\right), \quad 
h_c(u_{c}-c)=m, \quad m=-ch_{ \infty},
\end{equation}
and
\begin{equation} 
\sigma=m( h_c^2)'/6, \quad \varphi'=c+\frac{m}{h_c}\left(1+\frac{( h_c^2)''}{6}\right).
\end{equation}

\section{Numerical approach}
In this section we summarize the numerical approach to be applied 
in this paper. As in \cite{DK}, we will use a Fourier spectral 
method with the Krylov subspace technique GMRES \cite{GMRES}. 

The basic idea is to solve the SGN equations in Fourier space and to 
approximate the Fourier transform with a discrete Fourier transform 
(DFT)
which is computed with a fast Fourier transform (FFT). We use the 
following convention for the Fourier transform for sufficiently regular and localized functions $g$:
\begin{align*}
  &\forall (k_{x},k_{y})\in\mathbb{R}^{2}, \quad  
  \widehat{g}(k_{x},k_{y}) & :=\int_{\mathbb{R}^{2}}    
  e^{-i(k_{x}x+k_{y}y)}\,g(x,y) \, dxdy
    \\
 &\forall (x,y)\in\mathbb{R}^{2}, \quad    g(x) & 
    =\frac{1}{(2\pi)^{2}}\int_{\mathbb{R}^{2}}
    e^{i(k_{x}x+k_{y}y)}\,\widehat{g}(k_{x},k_{y}) \,dk_{x}dk_{y}.
\end{align*}
It is well known that smooth rapidly decreasing or smooth periodic 
functions $g(x,y)$
in each variable  can be expanded in a Fourier series in both $x$ 
and $y$, and that the coefficients of this series decrease rapidly in both 
indices. The DFT can be seen as a truncated Fourier series, and the 
numerical error in truncating is of the magnitude of the first 
neglected coefficients, see for instance \cite{trefethen} and references 
therein. In an abuse of notation, we will in the following denote the 
DFT of a function $g(x,y)$ ($x$ and $y$ being some vectors of 
collocation points) also by $\widehat{g}$. Derivatives are 
approximated in standard way, e.g., $\widehat{g_{x}}\approx 
ik_{x}\widehat{g}$. 

This approach is numerically efficient and of high accuracy for functions  vanishing  rapidly at infinity or being periodic. 
Therefore we will establish equations for $\mathbf{v}=\nabla \varphi$ 
having in contrast to $\varphi$ the property to vanish at infinity. 
For this, we replaced $\nabla\varphi$ by $\mathbf v$ excepting \eqref{2}  written in the  form 
\begin{equation}
\mathbf{v}_t+\nabla\left(\frac{\vert\mathbf{v}\vert^2}{2}-\frac{\vert \nabla\sigma\vert^2}{2h^2}+gh-\frac{9}{2}\frac{\sigma^2}{h^4}\right)=0.
\end{equation} 
If, initially, $\rm curl \mathbf{v} =0$, it will be vanishing for any 
$t>0$. Since derivatives in a Fourier spectral method are computed by 
multiplications with the dual variable in Fourier space, this is also 
the case for the numerical approach. On the other hand we are interested in minimizing the number of elliptic 
PDEs to be solved.  If we used $\eqref{eq:SGN-mom}$ for  $h\overline {\mathbf{u}}$ (note that the flux in $\eqref{eq:SGN-mom}$ 
depends of $\ddot h$), two such elliptic
PDEs would need to be solved for each time step instead of one for $\sigma$. Last not 
least we want to avoid the computation of high order derivatives as 
much as possible since this will lead to rounding errors that can 
pile up during the time evolution. \\
For these reasons we work with the variables $h$ and  $\mathbf{v}$ in 
the following and compute $\sigma$ in each time step via
\begin{equation}
	\frac{3\sigma}{h^{3}}-\nabla \cdot\left(\frac{\nabla \sigma}{h}\right) 
	= - \nabla \cdot \mathbf{v}
	\label{sigmav}.
\end{equation}
This is done as in \cite{DK} in Fourier space with the Krylov subspace technique GMRES 
\cite{GMRES}. 
Concretely we discretize in $x$ and $y$ in standard fashion for the 
DFT in each variable which means we introduce the collocation points 
$x_{n}=-L_{x}\pi+n 2\pi L_{x}$, $n=1,\ldots,N_{x}$ and 
$y_{m}=-L_{y}\pi+m 2\pi L_{y}$, $m=1,\ldots,N_{y}$. Thus we work on 
$\mathbb{T}^{2}$ with periods $2\pi L_{x}$ in $x$-direction and $2\pi 
L_{y}$ in $y$-direction. The positive numbers $L_{x}$, $L_{y}$ are 
chosen such that they correspond to the actual period in the 
respective direction in the case of a periodic solution or, for 
rapidly decreasing functions, that the function and its relevant 
derivatives vanish within the finite numerical precision. After this 
discretisation, the quantities $h$, $v_{x}$ and $v_{y}$ become 
$N_{x}\times N_{y}$ matrices. Equation (\ref{sigmav}) takes the form $\mathcal{L} 
\sigma=b$ where $b$ is a vector built from the matrix $\nabla \cdot 
\mathbf{v}$ after discretisation in a natural way by just writing the 
columns of the matrix in the form of a long vector. In the same way 
$\mathcal{L}$ is a $N_{x}N_{y}\times N_{x}N_{y}$ matrix the inverse 
of which is computed approximately via GMRES. Note that it is 
important to use a preconditioner $\mathcal{M}$ such that GMRES 
computes the solution of 
$\mathcal{M}^{-1}\mathcal{L}\sigma=\mathcal{M}^{-1}b$. We choose as a 
preconditioner $\mathcal{L}$ for $h=1$, the asymptotic value, i.e., 
$\mathcal{M}=3+k_{x}^{2}+k_{y}^{2}$. 

For the time integration of the SGN equations we apply as 
in \cite{DK} the standard explicit Runge-Kutta method of fourth 
order. 
This can be done directly for the first equation of (\ref{eq:SGN}). 
The equation  (\ref{2}) 
is written as $\varphi_{t}=F$ and then differentiated with respect to 
$x$ and $y$ by multiplying with $k_{x}$ and $k_{y}$ in Fourier space. 
This yields equations for $\mathbf{v}_{t}$ which are 
again integrated with the fourth order Runge-Kutta method. Note that 
GMRES is started with the last computed value for $\sigma$ as an 
initial guess. For small time steps, $\sigma$ does not change much 
over one time step, and consequently GMRES converges  rapidly. 

Note, however, that this is only the case as long as there is enough 
resolution in space as indicated by the decrease of the DFT 
coefficients. If they no longer decrease to machine precision during 
time evolution, GMRES will converge more slowly or not to the wanted 
accuracy, and rounding errors will pile up. To avoid such 
scenarios, we use a \textit{Krasny filter} \cite{Kra}, i.e., we put all 
DFT coefficients with modulus smaller than $10^{-12}$ equal to zero. 
However, it is recommended to provide enough resolution in space 
during the whole computation since GMRES will become slower in 
underresolved cases. Even if there is enough resolution, GMRES will take around 
90 \% of the computation time. 

\begin{figure}[!htb]
\includegraphics[width=0.49\hsize]{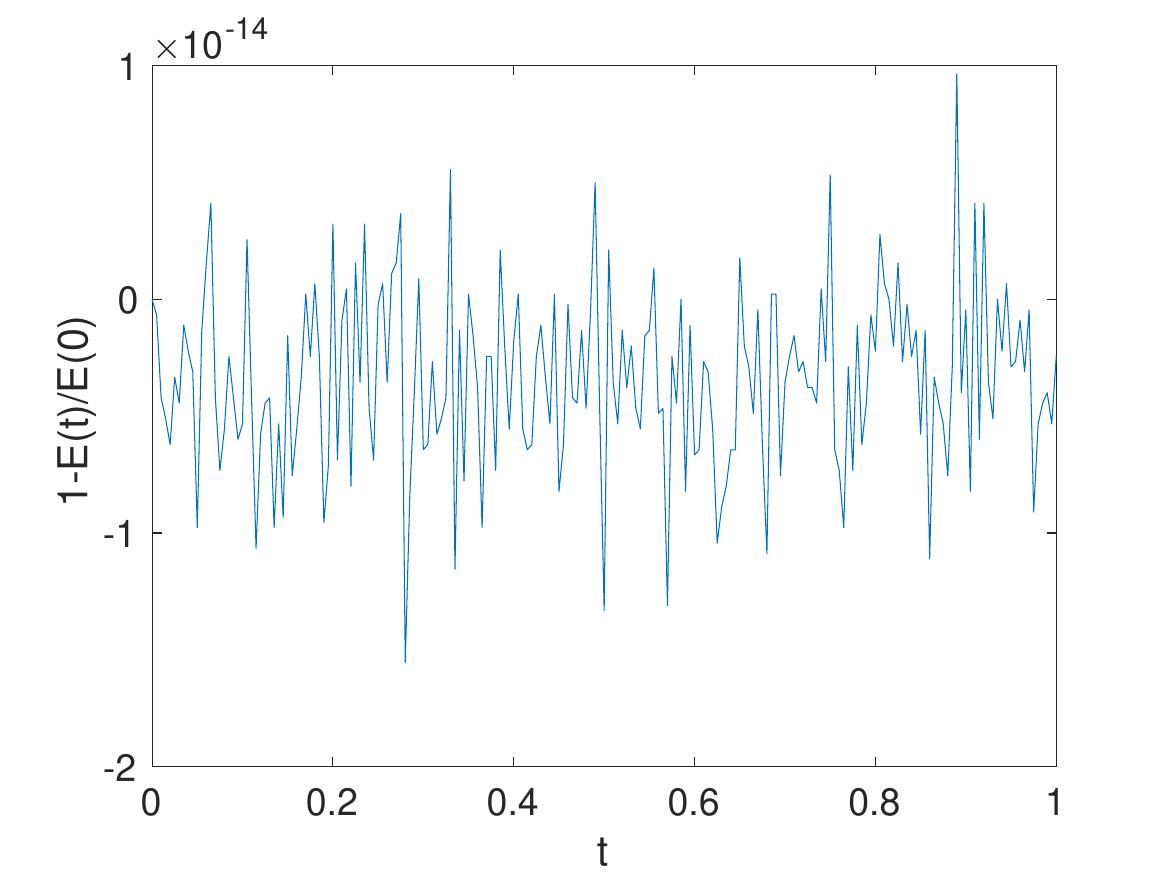}
\includegraphics[width=0.49\hsize]{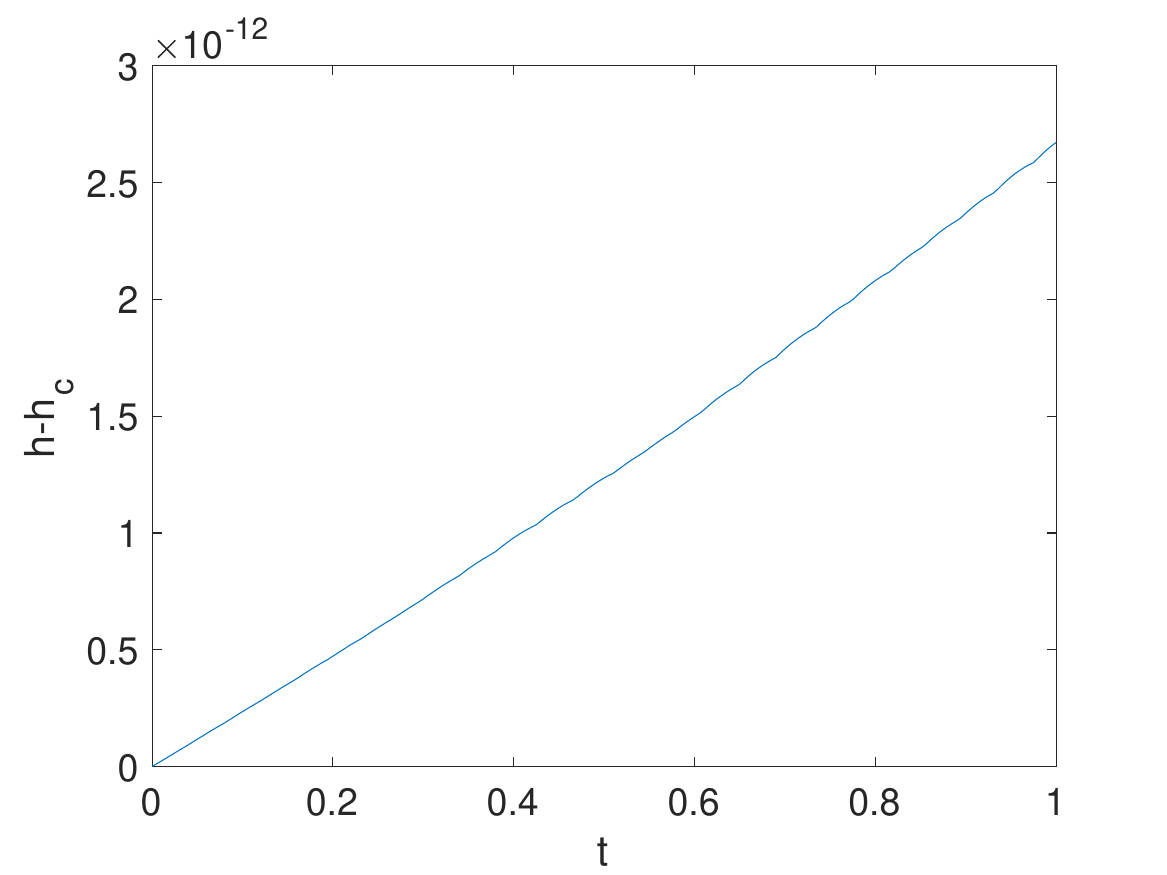}
\caption{Propagation of the line soliton with $c=1.7$.   We show on the left 
the relative conservation of the numerically computed energy,  and on the 
right the difference between exact and numerical solution for $h$. } 
\label{Delta}
\end{figure}

The 
accuracy of the time integration is controlled via the conservation 
of the numerically computed energy. Due to unavoidable numerical 
errors, the computed energy will depend on time, and the relative 
conservation of the energy can be used as discussed for instance in 
\cite{etna} to estimate the resolution in time: typically this 
quantity overestimates the numerical accuracy by 2 orders of 
magnitude. In all experiments reported below, relative energy 
conservation is better than $10^{-5}$ and thus above plotting 
accuracy. Since no explicit full 2D solution to the SGN equations in
the class of functions studied here is known, we can only test the 1D 
solutions. If we propagate the exact 1D soliton with  $N_{x}=2^{10}$ and $N_{y}=2^{7}$ Fourier modes with $L_{x}=10$ 
and $L_{y}=2$ and $N_{t}=10^{3}$ time steps for $t\leq 1$, we get 
energy conservation to the order of $10^{-14}$ as can be seen on the 
left of Fig.~\ref{Delta}, essentially the optimum that can be 
achieved with double precision. On the right of the same figure, we 
show the difference between numerical and exact solution for $h$  as a function  of time. It can be seen that numerical errors in the time 
integration pile up as usual and increase with time, but that the 
solitary wave can be  propagated with machine precision.

\section{Transverse stability of line solitary waves}
The 1D solitary wave (\ref{sol}) can be seen as a $y$-independent solution 
to the 2D SGN equations, i.e., a traveling wave solution infinitely 
extended in $y$-direction. If the solution (\ref{sol}) is considered in dependence 
of $y$ instead of $x$, it can be seen as an $x$-independent solution to the 2D SGN 
equations.  We call such solutions  {\it line solitary waves}. In this section 
we consider perturbations of such line solitary waves (for simplicity we 
consider an infinite extension in $y$). Numerically, the situation is 
studied  on $\mathbb{T}^{2}$ with period $2\pi L_{x}$ in 
$x$-direction and $2\pi L_{y}$ in $y$-direction. 
We consider in this section the case $c=1.7$. Cases with smaller 
and larger velocity have been studied as well, but the results are very similar 
to what is shown below and will thus not be presented. We will always put 
$h_{\infty}=g=1$ in the following We always show $\overline{\mathbf{u}}$ 
since it is more important in applications though the code 
is set up for $\mathbf{v}$. 

In this section we study various perturbations of the line solitary 
waves and show that they are numerically stable giving strong 
evidence to the first part of the main conjecture. 

\subsection{Deformed line solitary wave}

\begin{figure}[!htb]
\includegraphics[width=0.49\hsize]{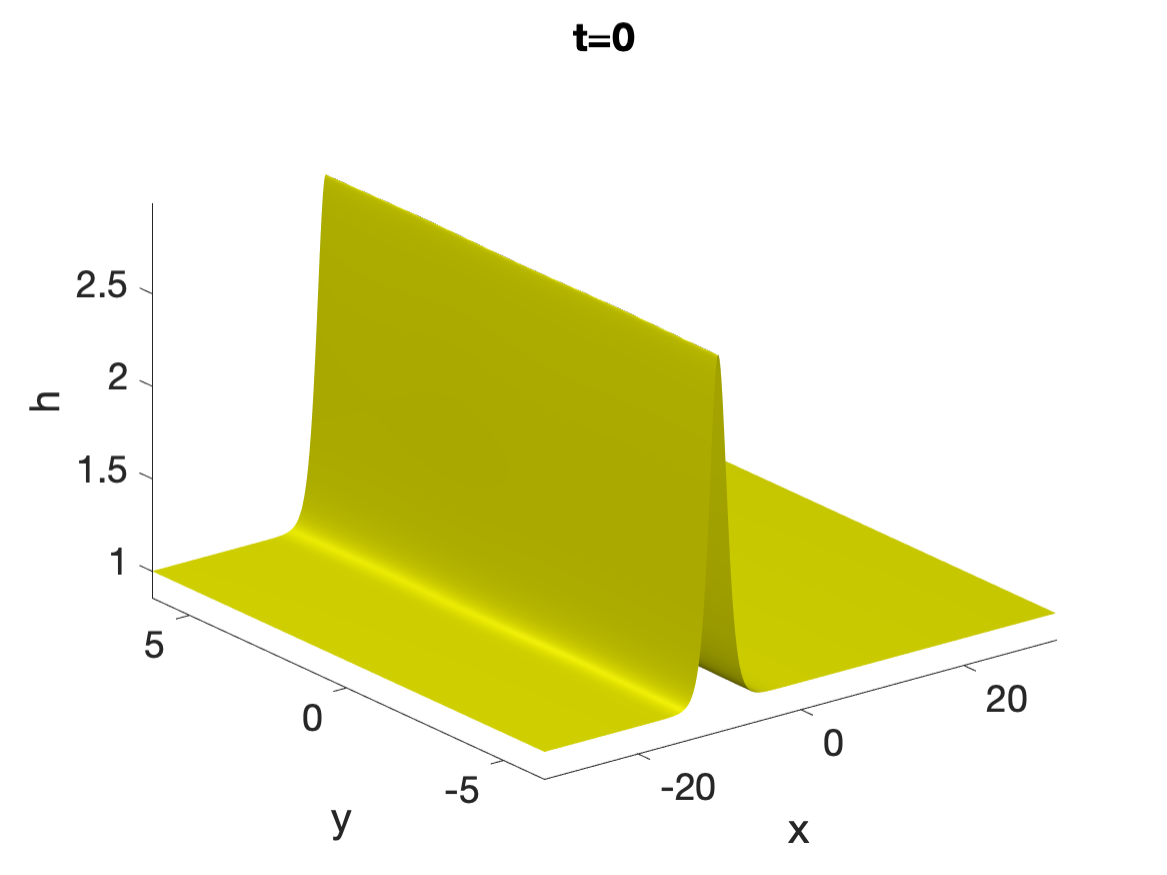}
\includegraphics[width=0.49\hsize]{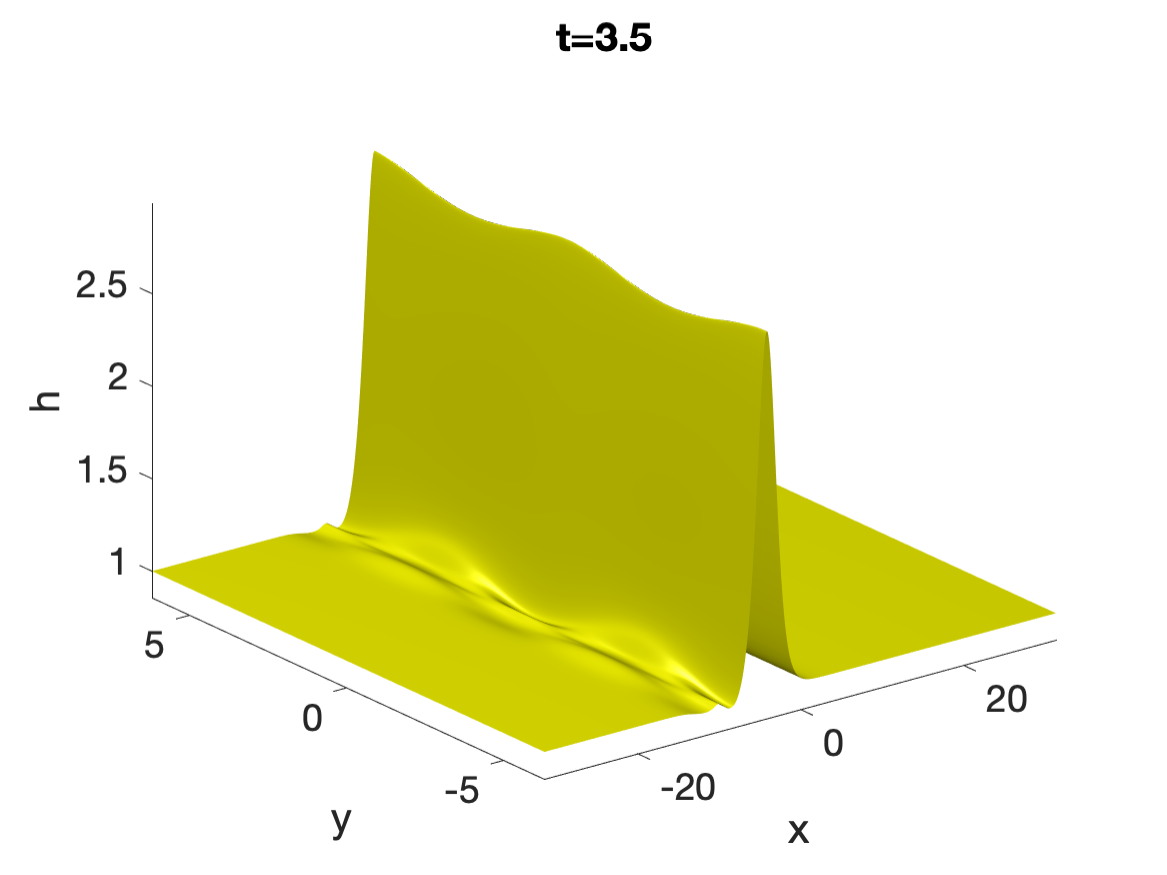}\\
\includegraphics[width=0.49\hsize]{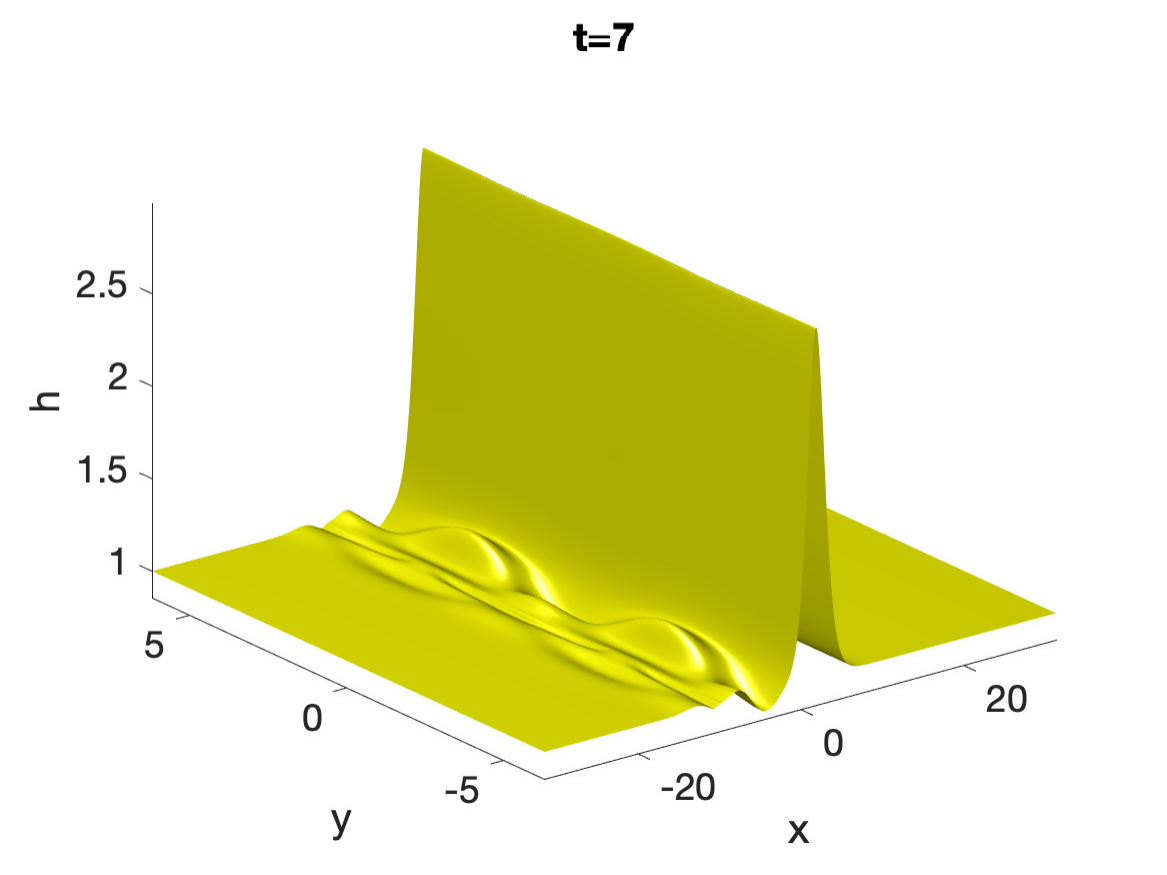}
\includegraphics[width=0.49\hsize]{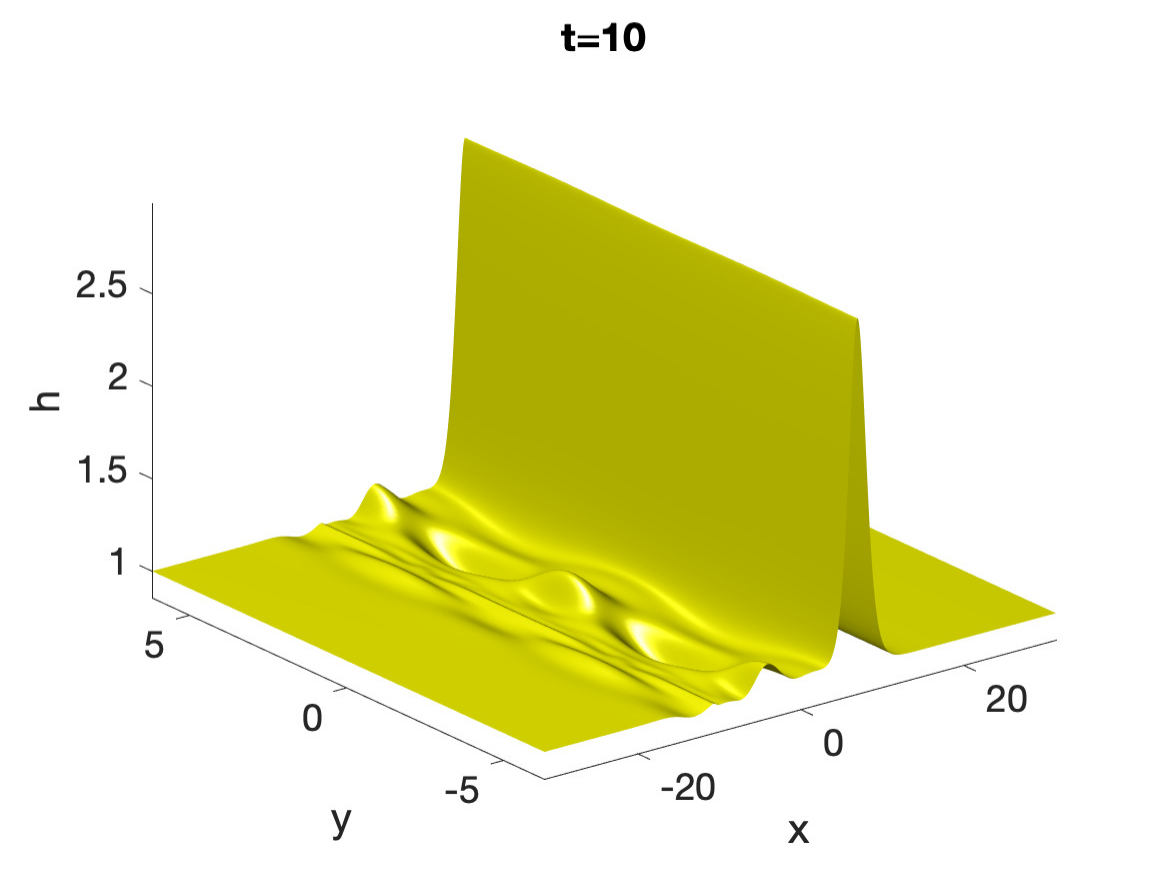}
\caption{Solution $h$ to the 2D SGN equation for initial data being a deformed line solitary 
wave of the form (\ref{soldef}) for several values of time. }
\label{soldefh}
\end{figure}

\begin{figure}[!htb]
\includegraphics[width=0.49\hsize]{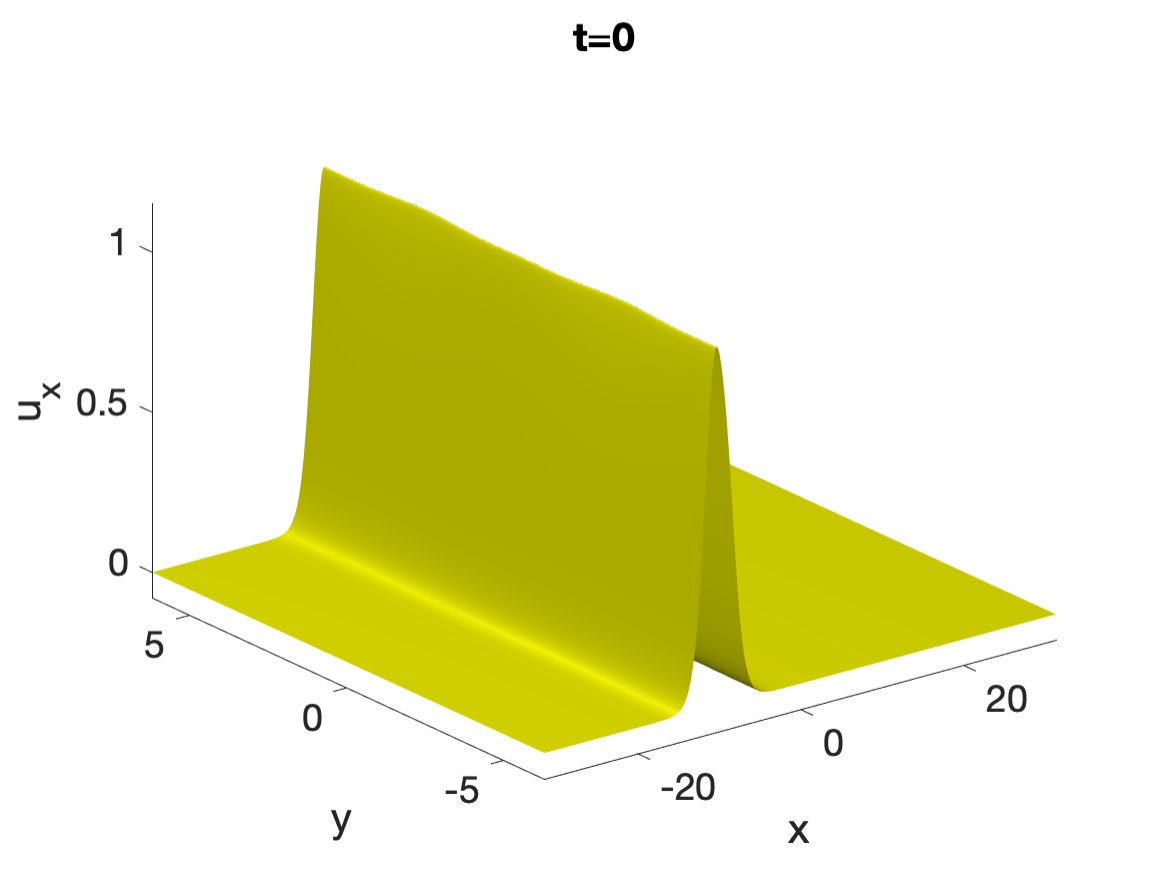}
\includegraphics[width=0.49\hsize]{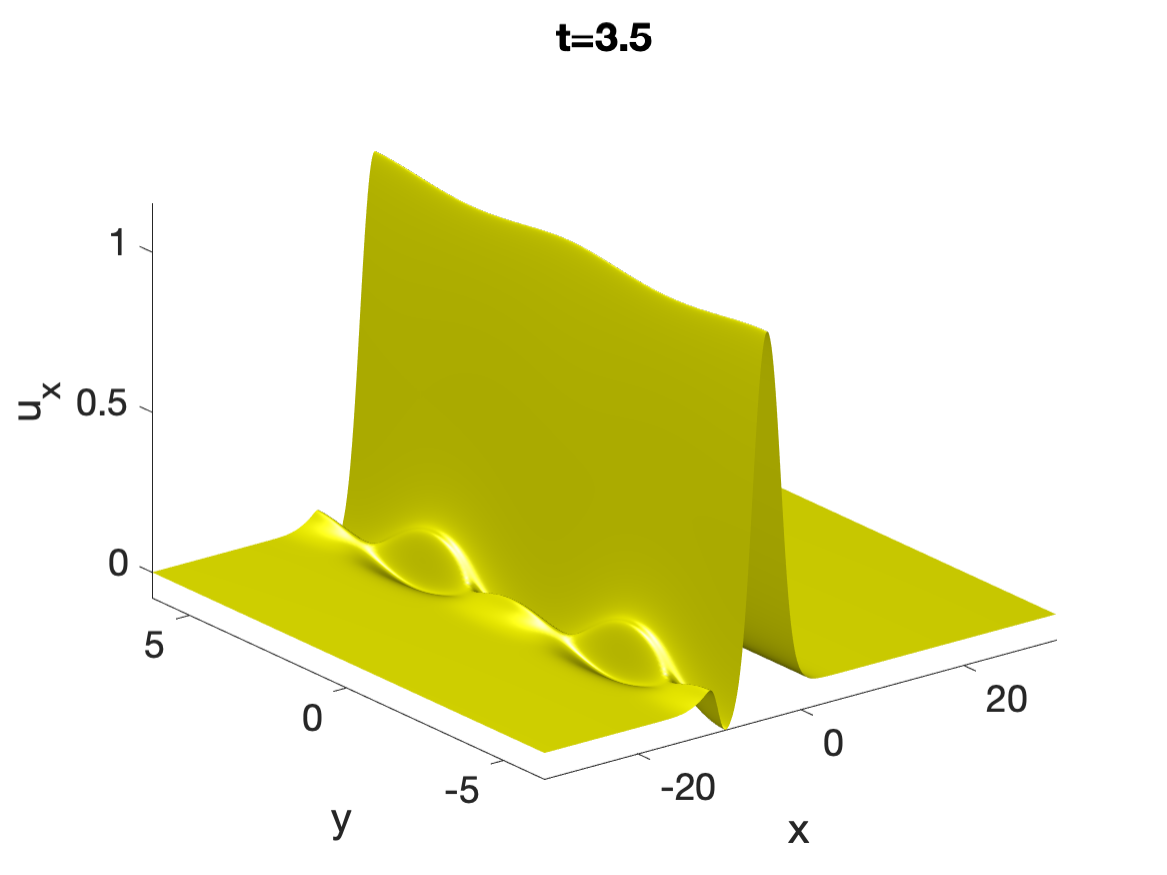}\\
\includegraphics[width=0.49\hsize]{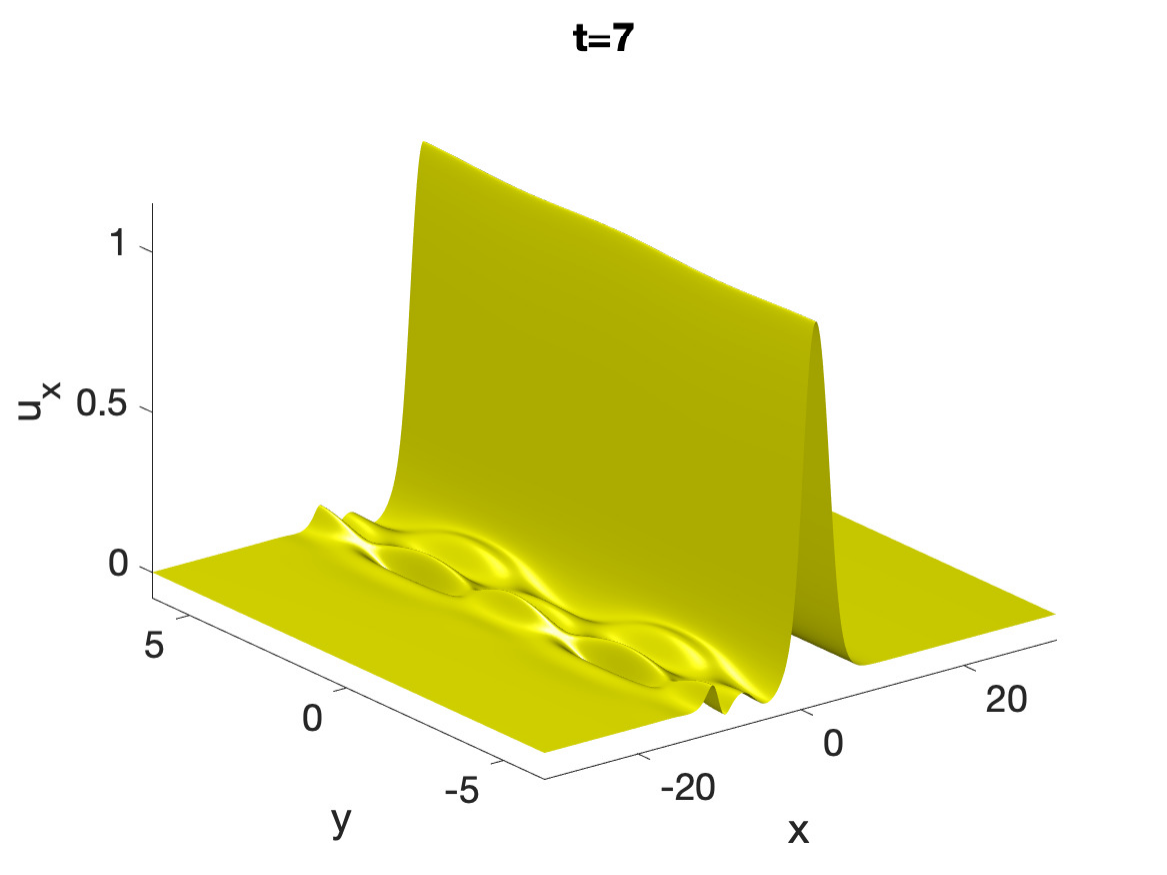}
\includegraphics[width=0.49\hsize]{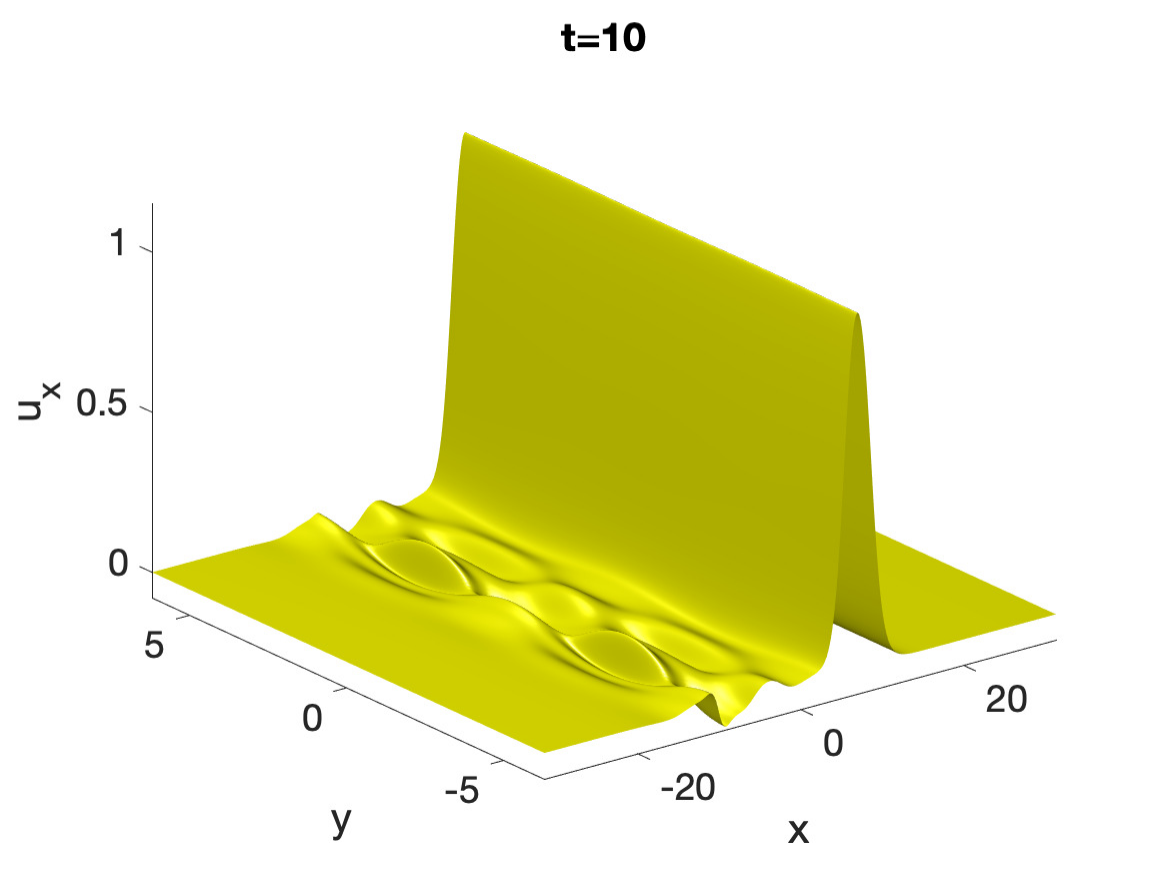}
\caption{Solution $u_{x}$ to the 2D SGN equation for initial data being a deformed line solitary 
wave of the form (\ref{soldef}) for several values of time. }
\label{soldefux}
\end{figure}

\begin{figure}[!htb]
\includegraphics[width=0.49\hsize]{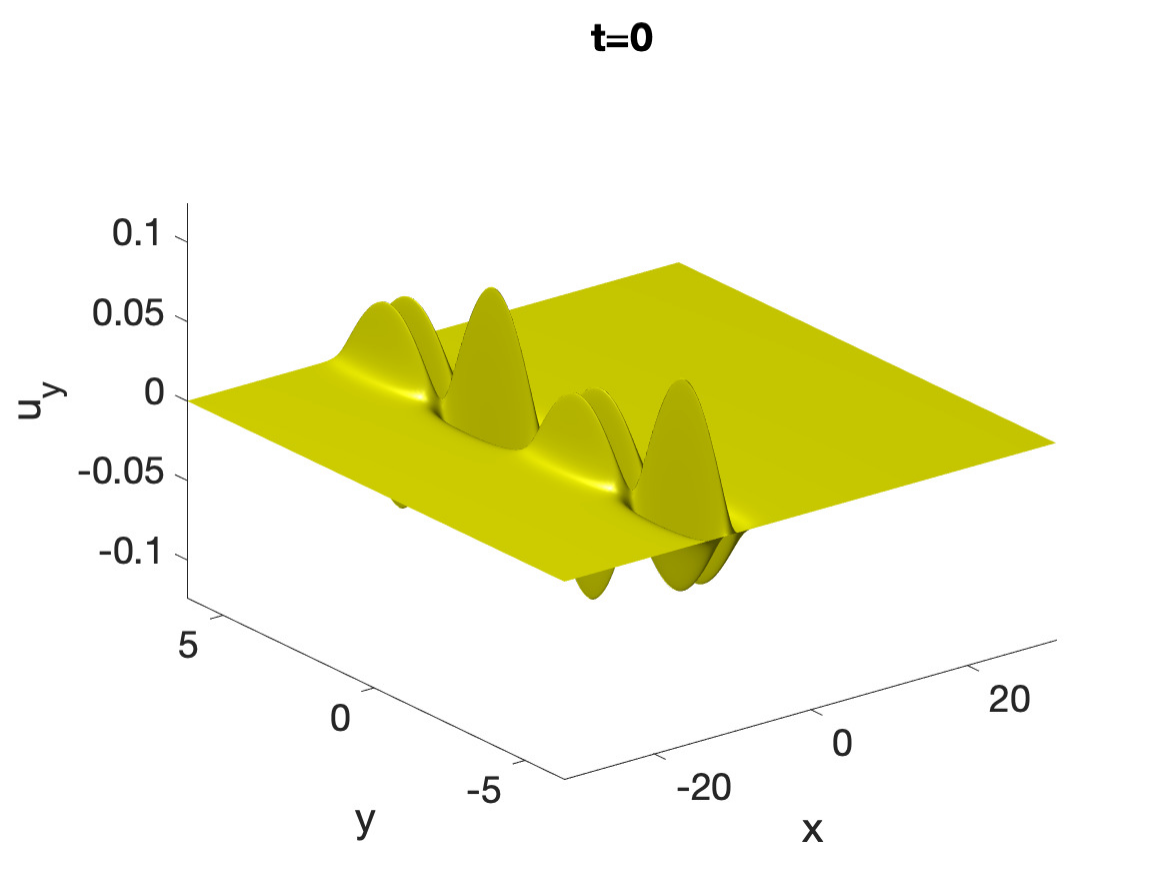}
\includegraphics[width=0.49\hsize]{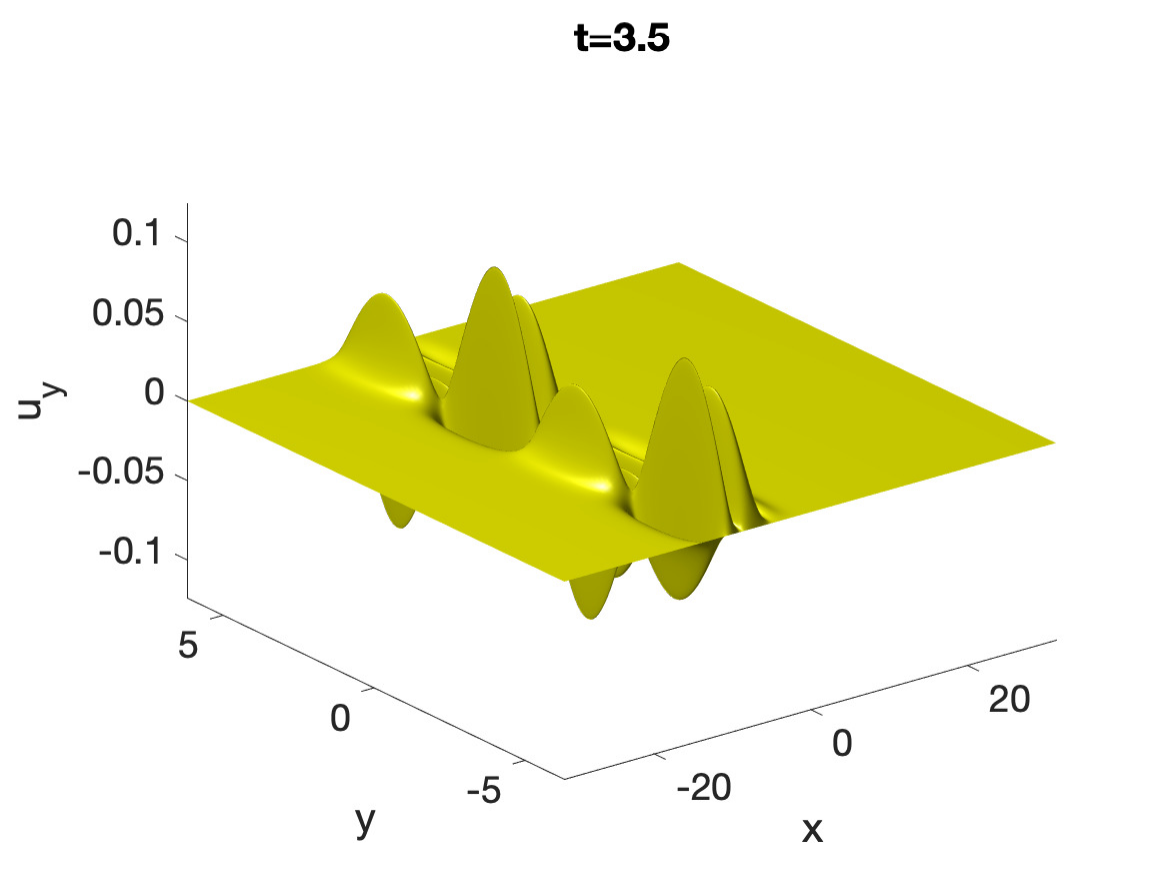}\\
\includegraphics[width=0.49\hsize]{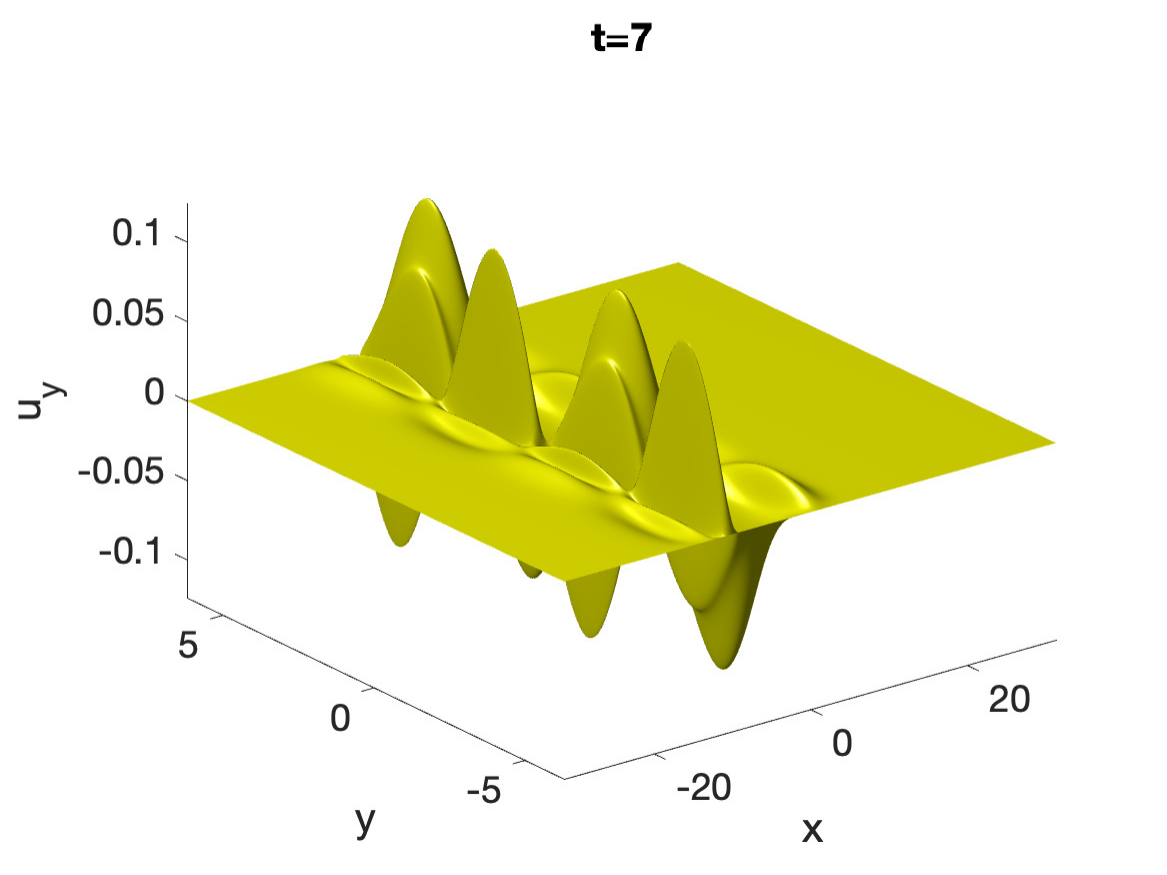}
\includegraphics[width=0.49\hsize]{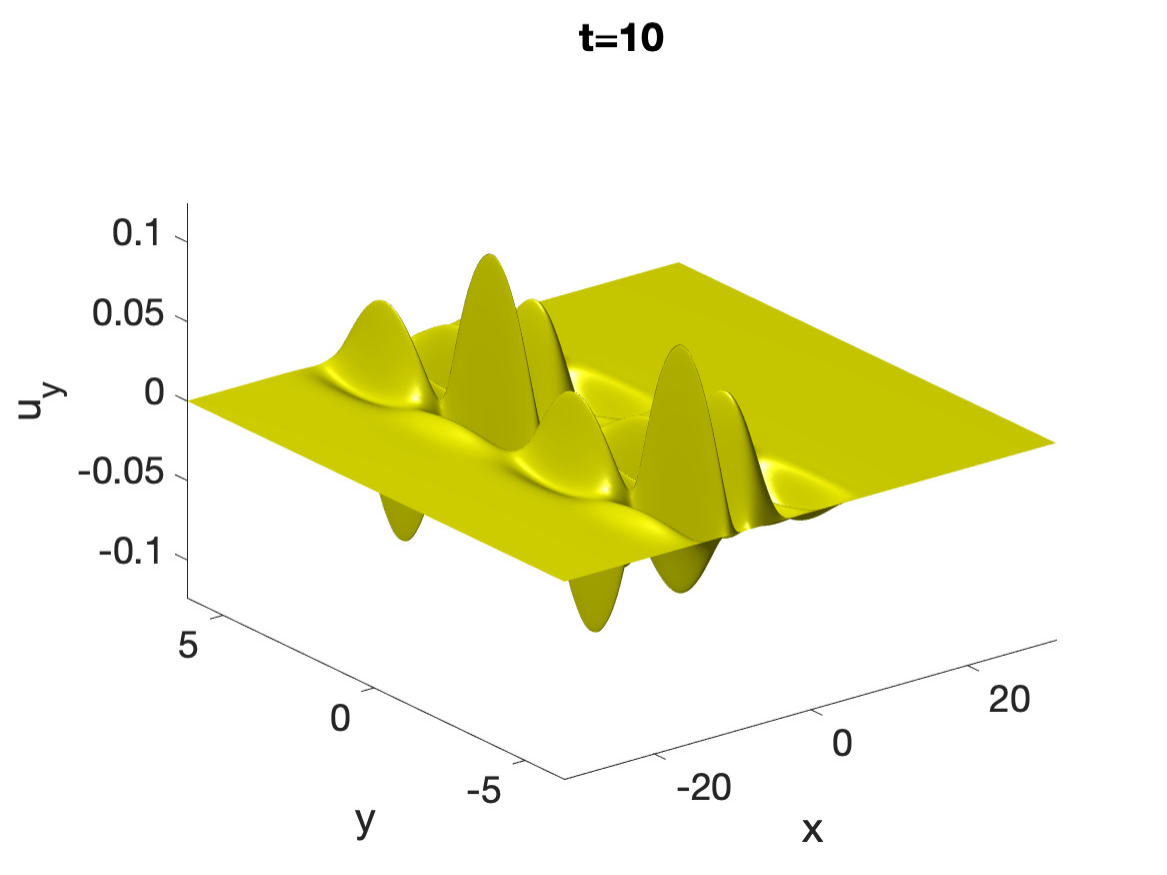}
\caption{Solution $u_{y}$ to the 2D SGN equation for initial data being a deformed line solitary 
wave of the form (\ref{soldef}) for several values of time. }
\label{soldefuy}
\end{figure}

As a first perturbation we consider a deformed line solitary wave,  i.e., 
initial data of the form 

\begin{subequations}
\label{soldef}
\begin{align}
h(x,y,0)&=1+(c^2-1)\mbox{sech}^2\left( \frac{x-x_{0}-0.1\cos(y)}{2c}\sqrt{3 (c^2-1)}\right), \\
 v_{x}(x,y,0)&=c-\frac{c}{h}\left(1+\frac{(h^2)_{xx}}{6}\right) ,\quad v_{y}(x,y,0)=0.
 \end{align}
\end{subequations}
Here $x_{0} =-10$.  The resulting initial data 
for $h$, $u_{x}$ and $u_{y}$ are shown as the first figure in 
Fig.~\ref{soldefh}, Fig.~\ref{soldefux} and Fig.~\ref{soldefuy} 
respectively. 

We use $N_{x}=2^{10}$ and $N_{y}=2^{7}$ Fourier modes with $L_{x}=10$ 
and $L_{y}=2$ and $N_{t}=10^{3}$ time steps for $t\leq 10$. The 
function $h$ for these initial data is shown in Fig.~\ref{soldefh} 
for several times. It can be seen that the initial deformation leads 
to a stronger one which eventually results in an emission of radiation 
from the deformed line solitary wave. The final state appears to be an 
unperturbed line solitary wave. 

The behavior of the quantity $u_{x}$ is very similar as can be seen 
in Fig.~\ref{soldefux}. The initial deformation leads to radiation 
plus the unperturbed  $u_{x}$ for the line soliton.

The quantity $u_{y}$ on the other hand vanishes in the unperturbed 
case and only shows in Fig.~\ref{soldefuy} radiative behavior during time evolution (the 
amplitude is always considerably smaller than the ones of $h$ in 
Fig.~\ref{soldefh} and $u_{x}$ in Fig.~\ref{soldefux}). 

\begin{figure}[!htb]
\includegraphics[width=0.49\hsize]{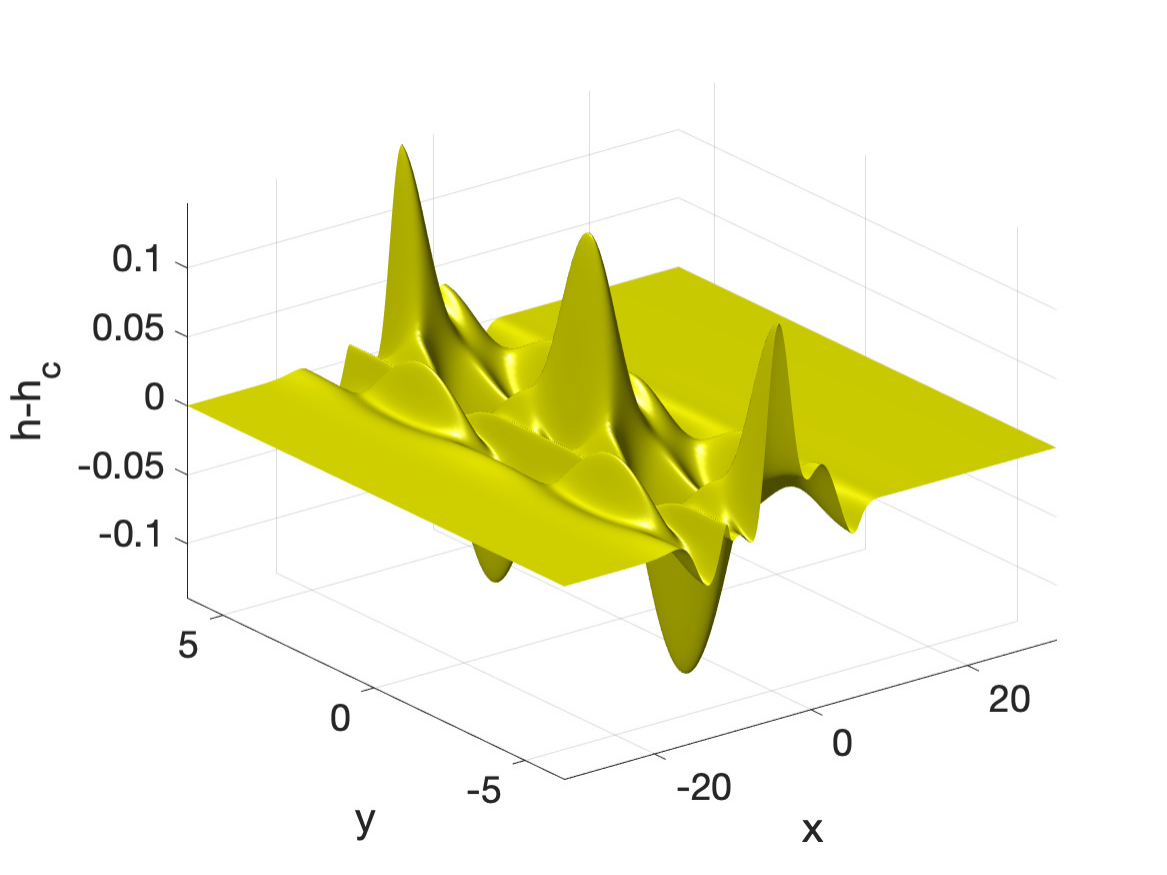}
\includegraphics[width=0.49\hsize]{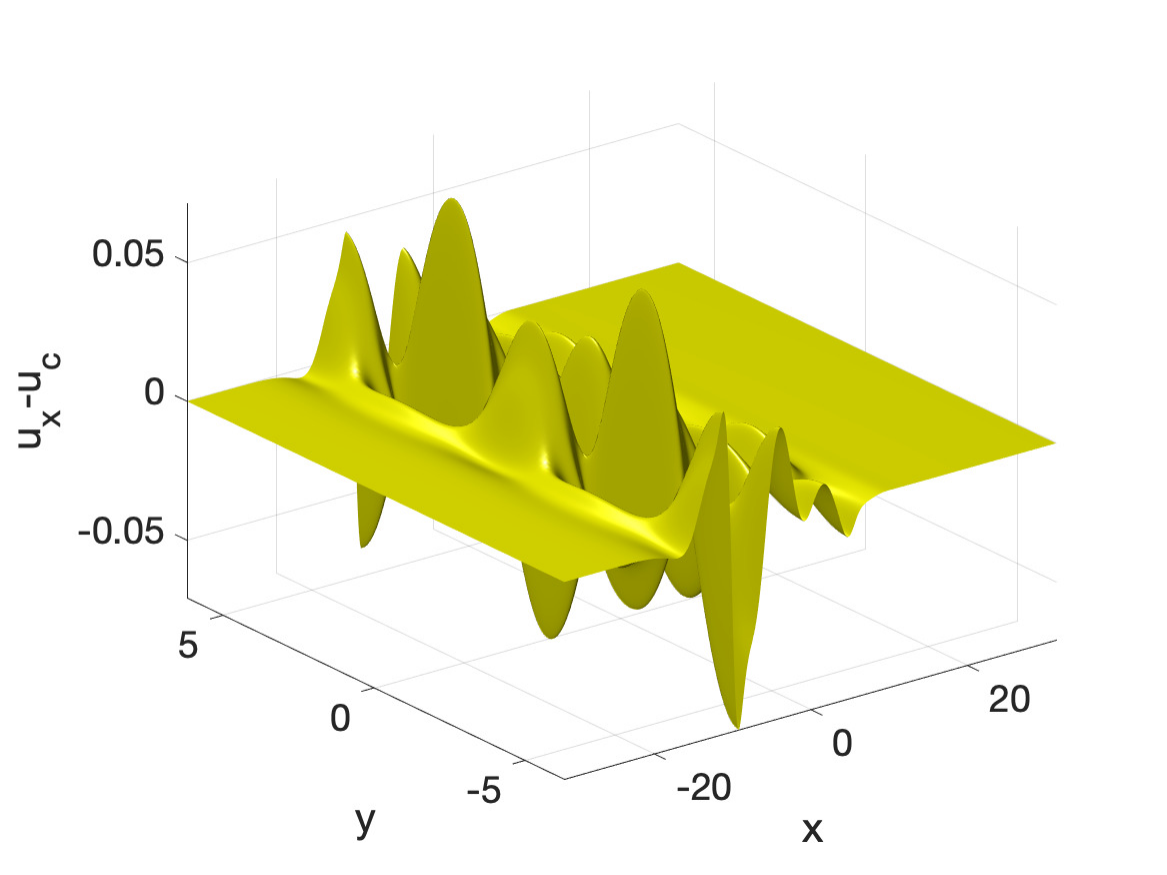}
\caption{Difference between the solution  to the 2D SGN equation for initial data being a deformed line solitary 
wave of the form (\ref{soldef}) for $t=10$ and an unperturbed line 
solitary wave centered at $x_{s}=6.995$: on the left for $h$, on the 
right for $u_{x}$. }
\label{soldefdiff}
\end{figure}

Note that the initial data (\ref{soldef}) are a perturbation of the 
order of $10\%$ and thus by no means small. However, it appears that 
the final state is the unperturbed line solitary wave which would 
indicate its stability. To make this more evident, we show in 
Fig.~\ref{soldefdiff} the difference between the final state and the 
unperturbed solitary wave (\ref{sol}) centered at $x_{s}=6.995$, the 
numerically determined location of the maximum of $h$. In 
Fig.~\ref{soldefdiff} we consider the 
difference between $h$ and the value $h_{c}$ for the solitary wave 
(\ref{sol})
with $c=1.7$ as well as the respective difference between $u_{x}$ and 
the unperturbed $u_{c}$ (\ref{sol}). It can be seen that this difference is of 
the order of the quantity $u_{y}$, no trace of the line solitary wave 
can be seen in the difference. This provides strong numerical evidence 
for the conjecture that a perturbation of the line solitary wave 
leads to the unperturbed solitary wave plus radiation for 
sufficiently large times. 

\subsection{Line solitary wave plus Gaussian perturbation}

We consider initial data with a localized perturbation of the form 

	\begin{subequations}
\begin{align}
	h(x,y,0) &= h_{c}(x-x_{0})\pm 0.1 \exp(-(x-x_{0})^{2}-y^{2}),\\
	v_{x}(x,y,0)&=c-\frac{
	c}{h}\left(1+\frac{(h^2)_{xx}}{6}\right),
	\quad v_{y}(x,y,0)=0,
\end{align}
	\label{solgauss}
\end{subequations}
i.e., a solitary wave with a small Gaussian perturbation. We work with 
$N_{x}=2^{11}$, $N_{y}=2^{8}$, $L_{x}=20$, $L_{y}=2$, $x_{0}=-20$ and 
$N_{t}=2000$ time steps for $t\leq 20$ (for larger times the same 
time step is used). The DFT coefficients 
decrease during the whole computation to the order of $10^{-7}$, 
relative energy conservation is of the order of $10^{-9}$. 

\begin{figure}[!htb]
\includegraphics[width=0.49\hsize]{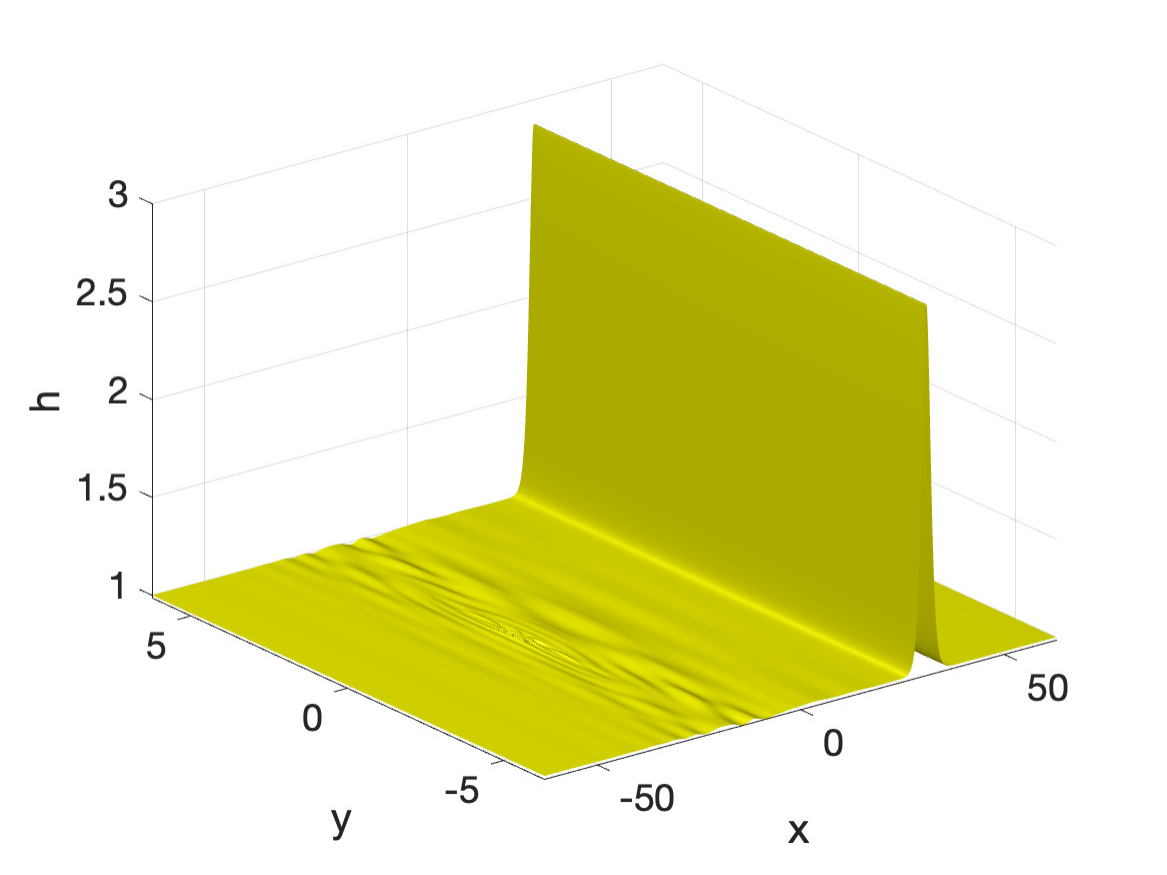}
\includegraphics[width=0.49\hsize]{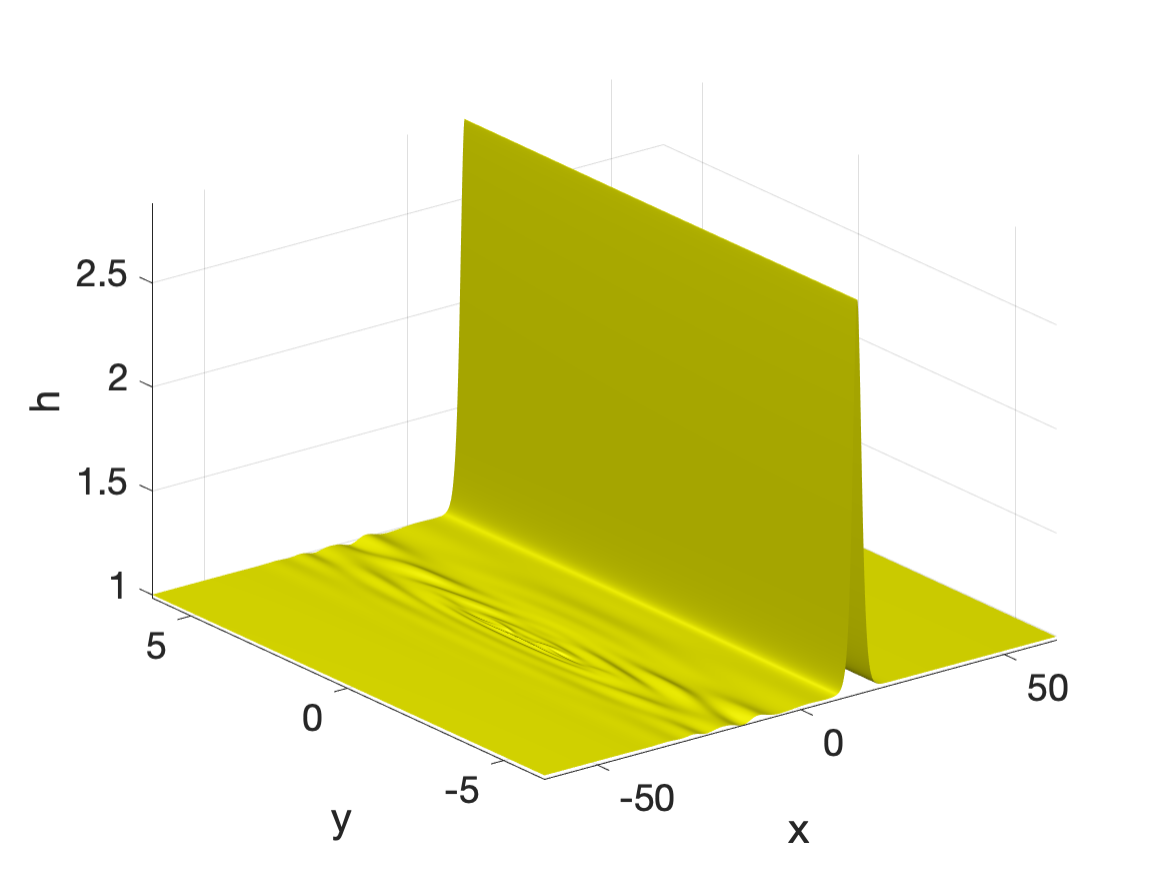}
\caption{Solution  $h$ to the 2D SGN equation for initial data of the 
form (\ref{solgauss}): on the left for the minus sign in 
(\ref{solgauss})  for $t=30$, on the 
right for the plus sign  for $t=20$. }
\label{solgaussh}
\end{figure}

The solution $h$ for both initial data (\ref{solgauss}) at the final 
time  can be seen in Fig.~\ref{solgaussh}. It appears that the 
final state is once more the unperturbed line solitary wave plus 
radiation. 

\begin{figure}[!htb]
\includegraphics[width=0.49\hsize]{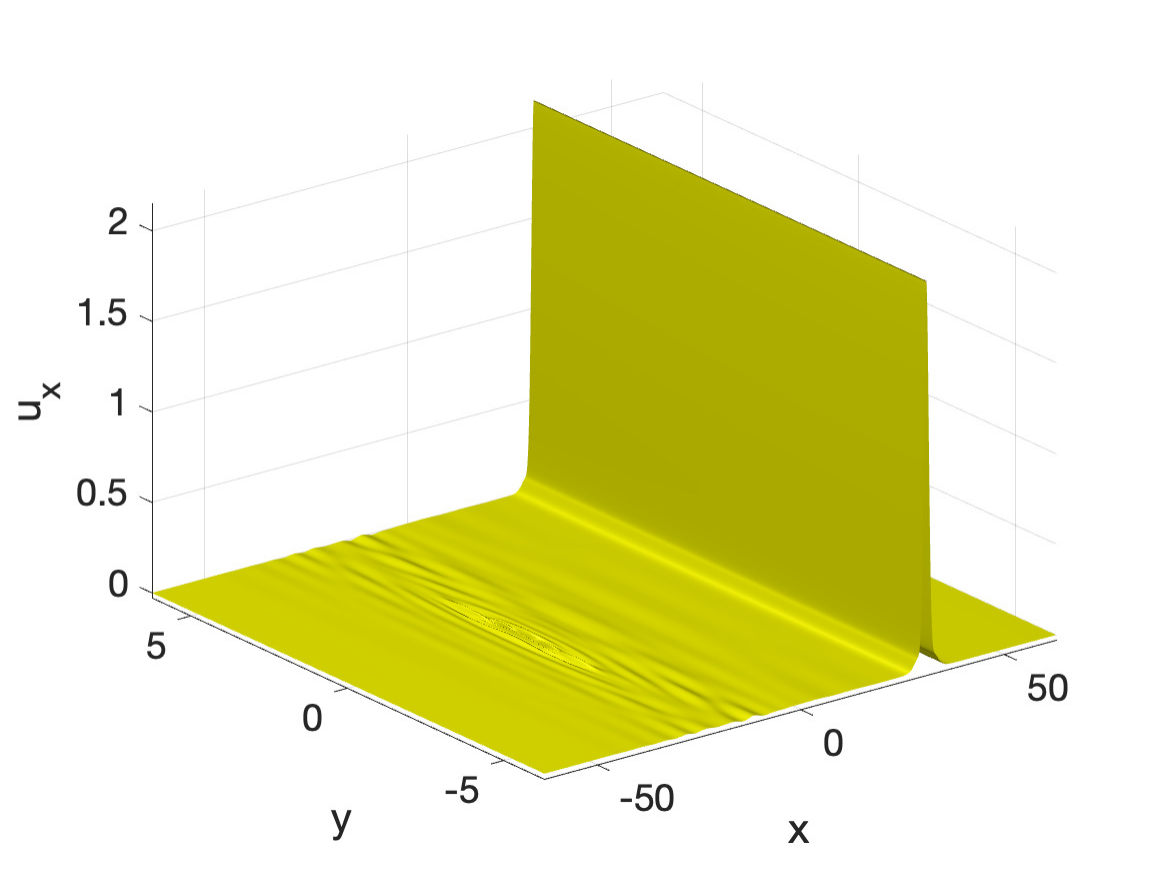}
\includegraphics[width=0.49\hsize]{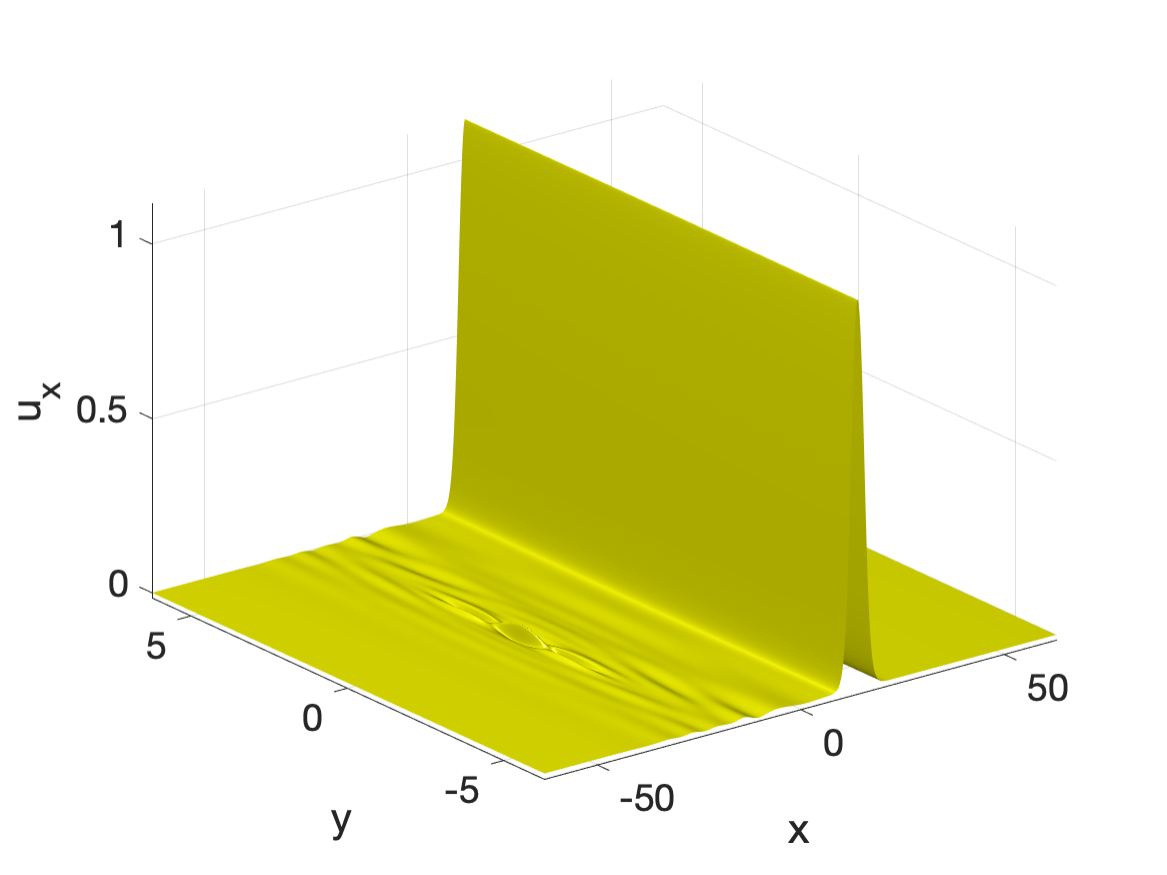}
\caption{Solution  $u_{x}$ to the 2D SGN equation for initial data of the 
form (\ref{solgauss}): on the left for the $-$ sign in 
(\ref{solgauss})  for $t=30$, on the 
right for the $+$ sign  for $t=20$. }
\label{solgaussux}
\end{figure}

The corresponding velocities $u_{x}$ at the final time are shown in 
Fig.~\ref{solgaussux}. Again this can be interpreted as a line 
soliton plus radiation. 

\begin{figure}[!htb]
\includegraphics[width=0.49\hsize]{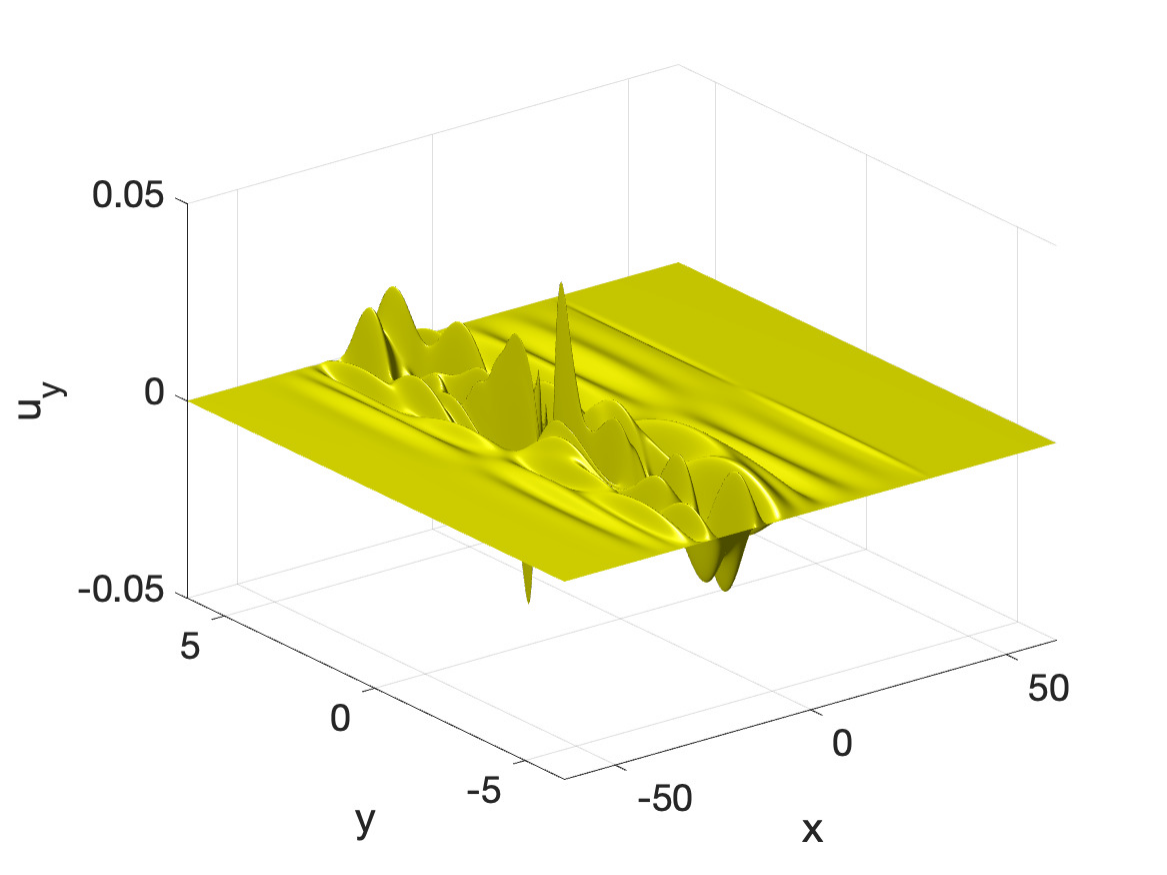}
\includegraphics[width=0.49\hsize]{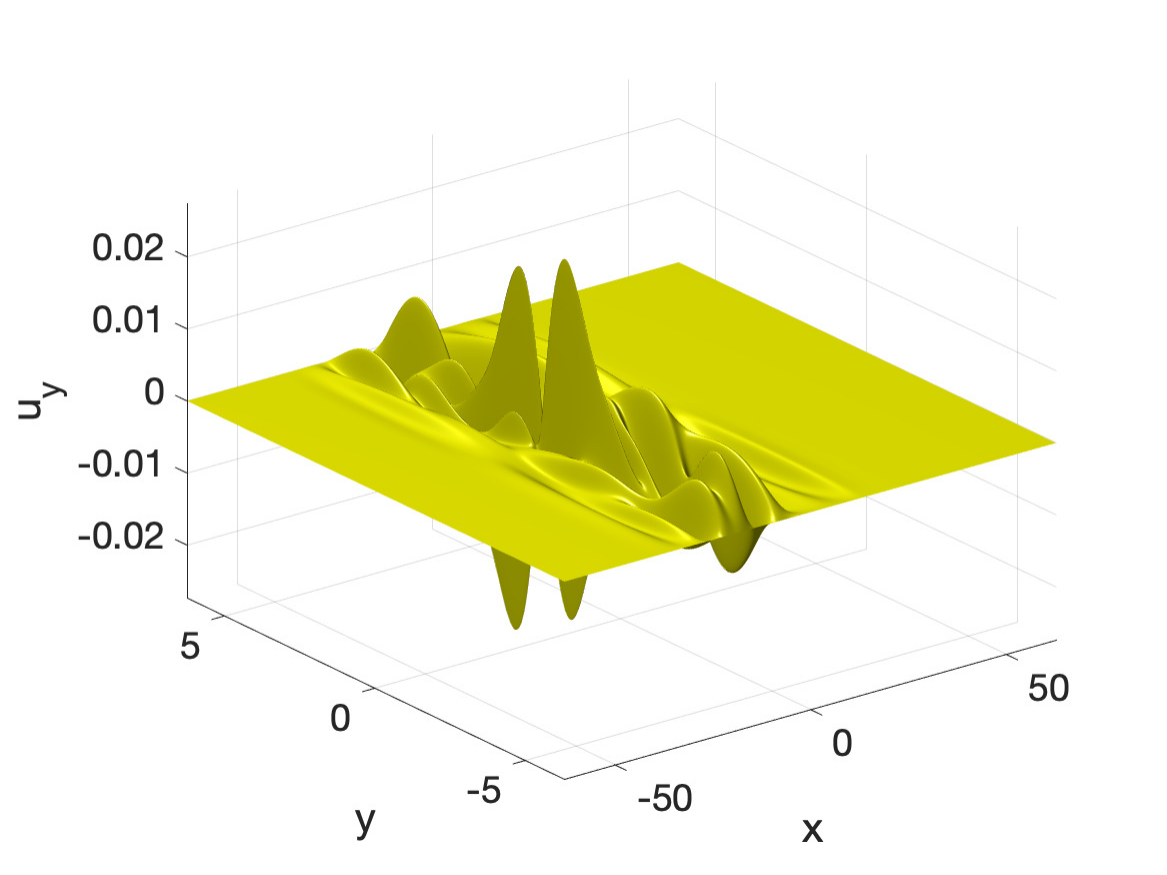}
\caption{Solution  $u_{y}$ to the 2D SGN equation for initial data of the 
form (\ref{solgauss}): on the left for the $-$ sign in 
(\ref{solgauss})  for $t=30$, on the 
right for the $+$ sign  for $t=20$. }
\label{solgaussuy}
\end{figure}

The velocities $u_{y}$ also show the expected behavior, i.e., only 
radiation, as can be seen in Fig.~\ref{solgaussuy}. Its magnitude is 
considerably smaller than the initial perturbation which stresses the 
stability of the line solitary wave. 

\begin{figure}[!htb]
\includegraphics[width=0.49\hsize]{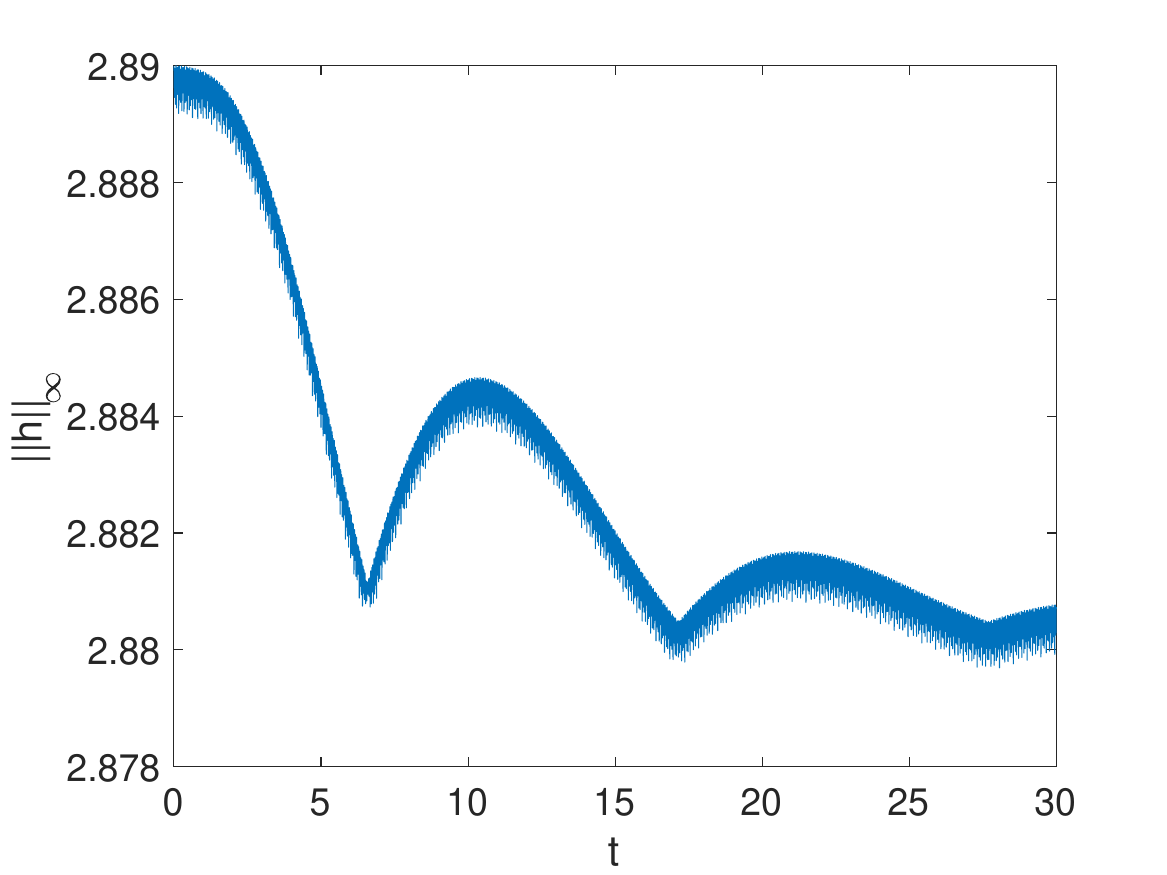}
\includegraphics[width=0.49\hsize]{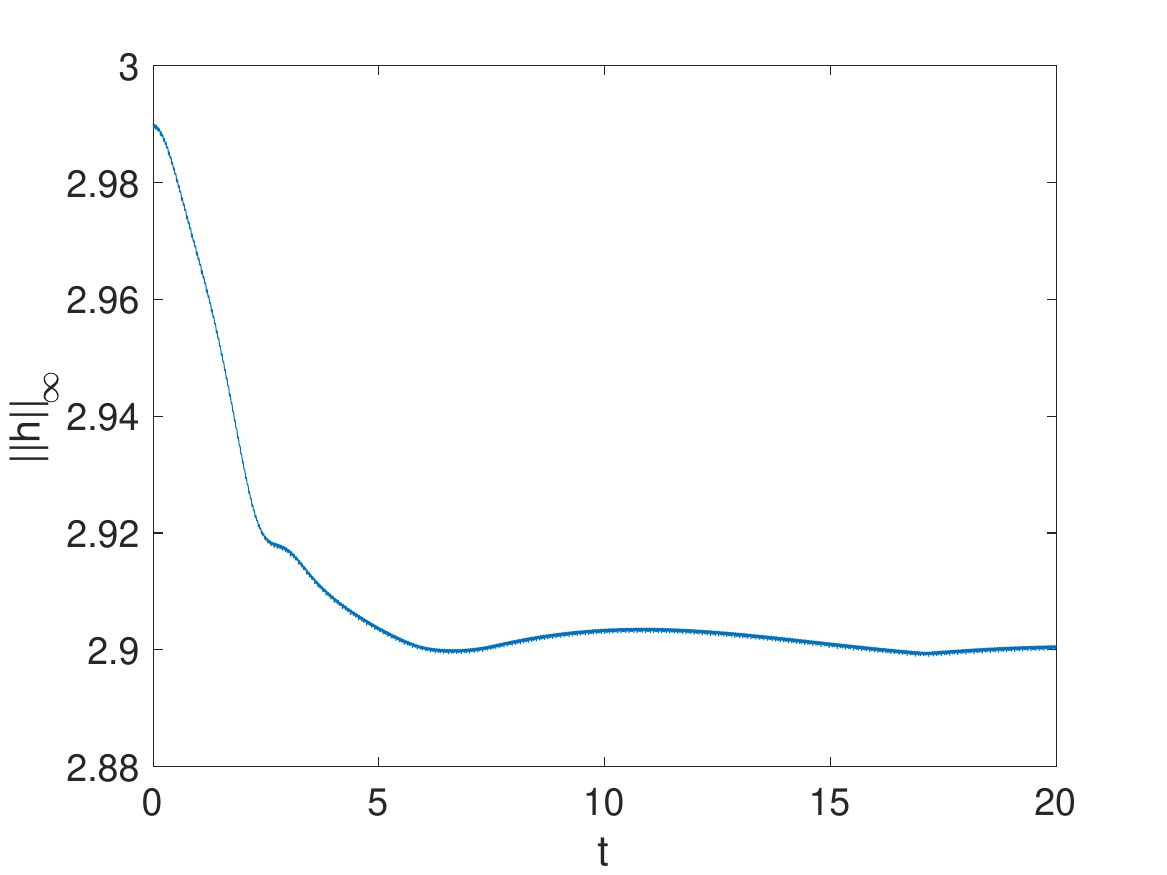}
\caption{$L^{\infty}$ norm of the solution  $h$ to the 2D SGN equation for initial data of the 
form (\ref{solgauss}) in dependence of time: on the left for the $-$ sign in 
(\ref{solgauss}), on the 
right for the $+$ sign. }
\label{solgaussmax}
\end{figure}

The $L^{\infty}$ norms of the solution $h$ can be seen in 
Fig.~\ref{solgaussmax}. They show the typical oscillations around a 
presumed final state (the small oscillations are due to the fact that 
we determine the $L^{\infty}$ norm on grid points which not 
necessarily correspond to the location of the maximum). Since we consider a perturbation with a higher 
respectively lower mass of several percents, the $L^{\infty}$ norm of 
this final state will have a slightly different value than the 
unperturbed soliton, and thus a slighly different velocity $c$. Note 
that there are stronger oscillations in the case with smaller mass. 
This indicates that one does not get as rapidly close to the final state as 
in the case with larger mass though the final time is the same, 
which is why we considered longer times in the case 
of a $-$ sign in (\ref{solgauss}). The 
final time is 
also at least twice than what was taken for the initial data (\ref{soldef}). 

\begin{figure}[!htb]
\includegraphics[width=0.49\hsize]{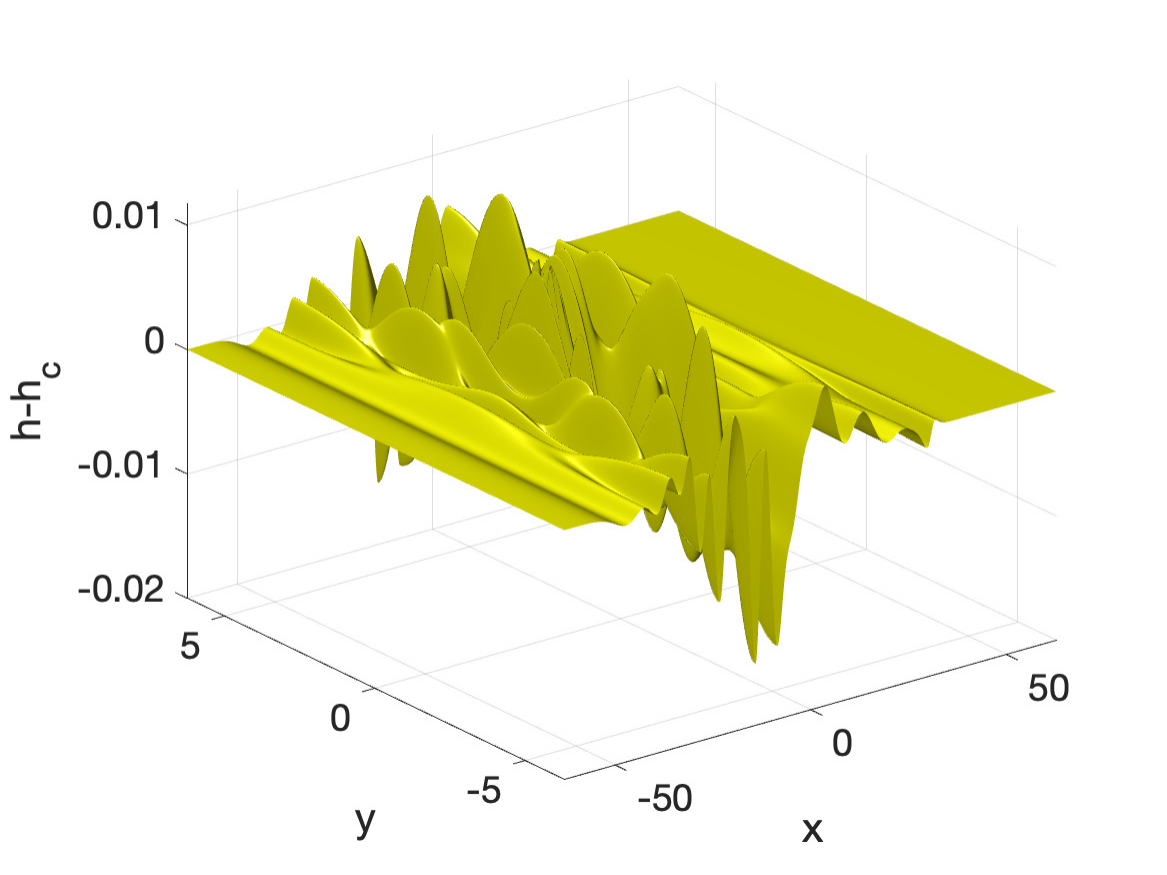}
\includegraphics[width=0.49\hsize]{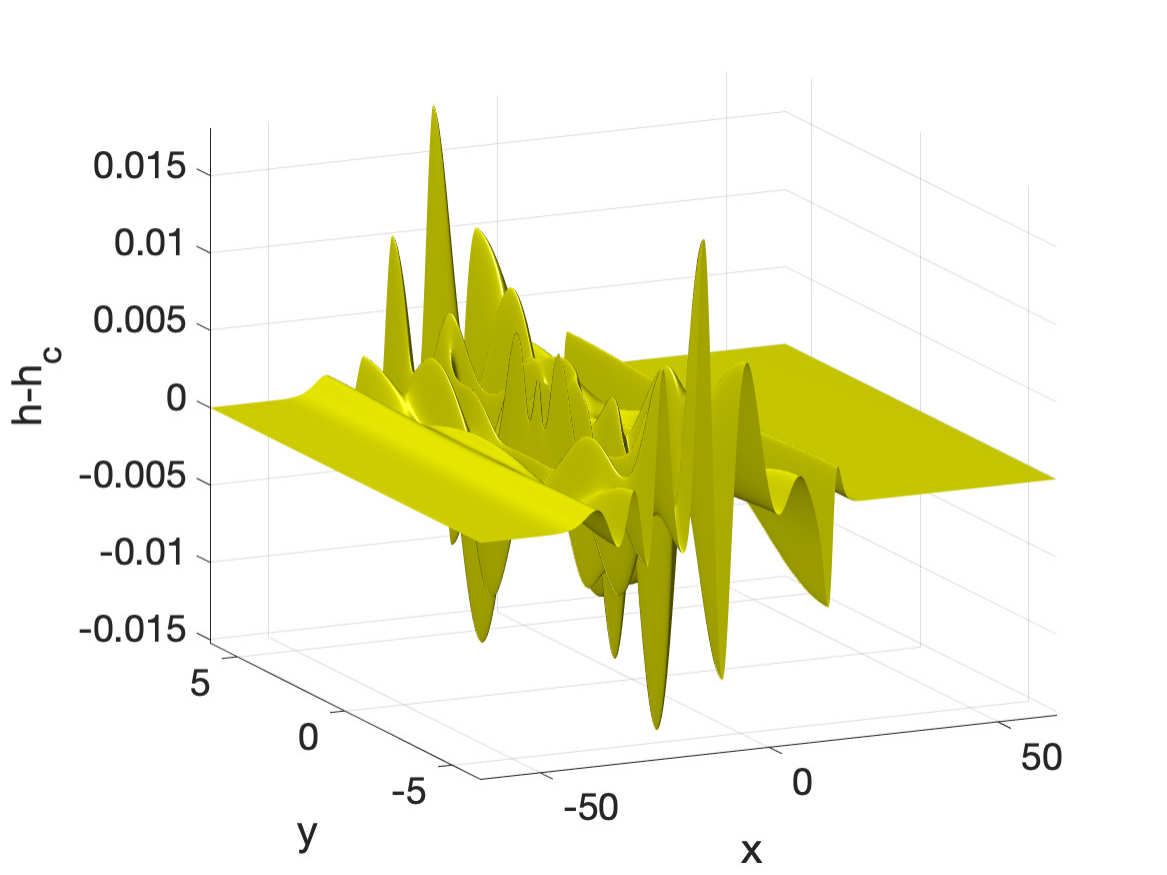}\\
\includegraphics[width=0.49\hsize]{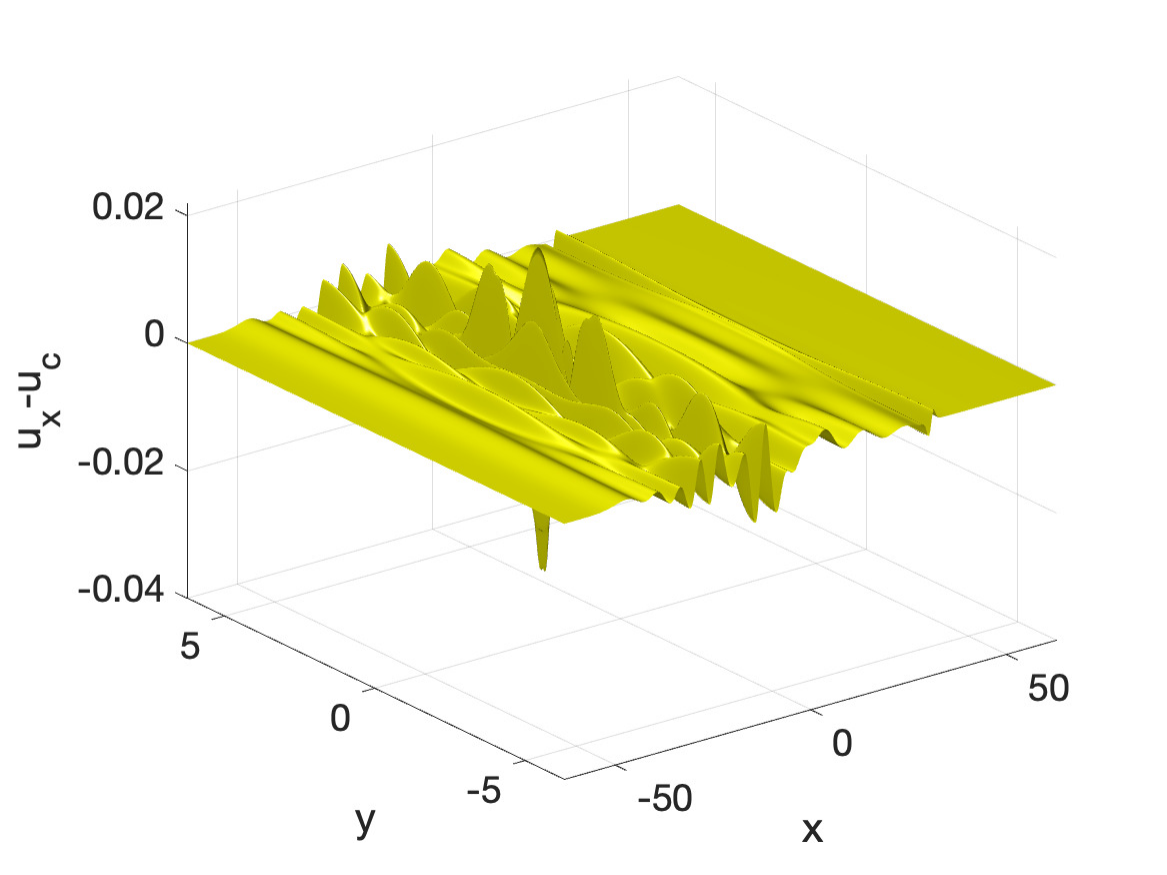}
\includegraphics[width=0.49\hsize]{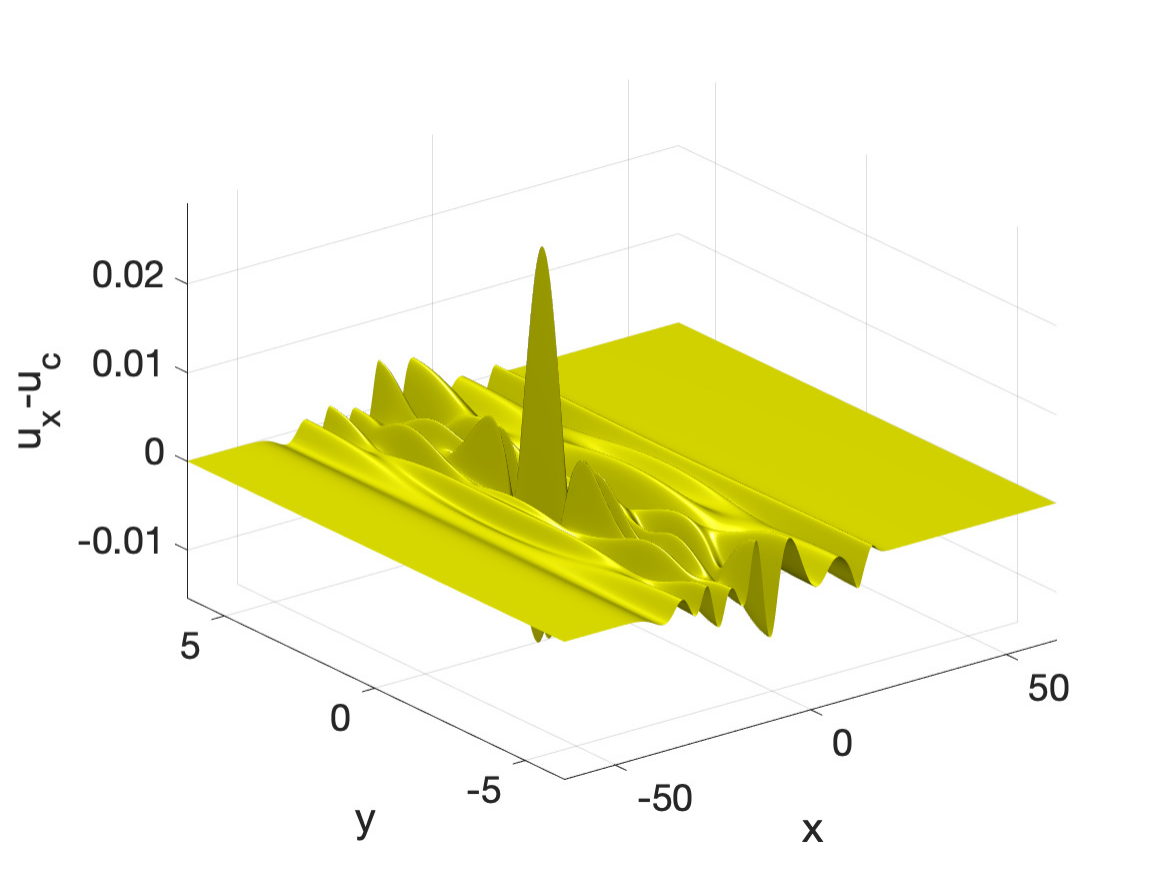}
\caption{Difference of the solution to the 2D SGN equation for initial data of the 
form (\ref{solgauss}) and a fitted line solitary wave: on the left for the $-$ sign in 
(\ref{solgauss}) for $t=30$, on the 
right for the $+$ sign  for $t=20$, the upper row for $h$, the lower for $u_{x}$. }
\label{solgaussdiff}
\end{figure}

In order to check whether the final state is indeed a line solitary 
wave, we consider once more the difference between the 1D solution 
(\ref{sol}) and the above 2D solution for $t=20$. The former is taken at 
the location of the maximum as before, the value of $c$ is obtained 
by fitting this maximum to the maximum of expression (\ref{sol}). The resulting differences for 
$h$ and $u_{x}$ are shown in Fig.~\ref{solgaussdiff}. It can be seen 
that the difference is of the order of magnitude of the radiation. As indicated by 
the $L^{\infty}$ norm in Fig.~\ref{solgaussmax},  the final state 
will be reached in the $-$ case in (\ref{solgauss}) only for even larger 
times.

\section{Initial data of crossing solitary waves}
Integrable equations in 2D as the KP II equation have 
multi-line-solitons as exact solutions, i.e. several line solitons 
crossing at certain angles forming time dependent patterns of 
solitons, see e.g. \cite{Kod} and references therein. Since SGN does 
not appear to be completely integrable, no exact multi-solitary waves 
are expected. But we can discuss the time evolution of initial data 
which are simply the superposition of one line soliton in 
$x$-direction and one in $y$-direction. Concretely we consider the 
following example,
\begin{subequations}
\begin{align}
h(x,y,0) & = h_{c}(x)+h_{c}(y)\\
v_{x}(x,y,0)&=c-\frac{ c}{h_{c}(x)}\left(1+\frac{((h_{c}(x))^2)_{xx}}{6}\right),\\
v_{y}(x,y,0)&=c-\frac{ c}{h_{c}(y)}\left(1+\frac{((h_{c}(y))^2)_{yy}}{6}\right).
\end{align}
\label{cross}
\end{subequations}
We consider the example $c=1.7$ for both line solitary waves which 
leads to a symmetric situation. We use $N_{x}=N_{y}=2^{10}$ Fourier 
modes with $L_{x}=L_{y}=10$ and $N_{t}=10^{3}$ time steps for $t\leq 
10$. The energy is conserved during the whole computation relatively 
to the order of $10^{-9}$. 

\begin{figure}[!htb]
\includegraphics[width=0.49\hsize]{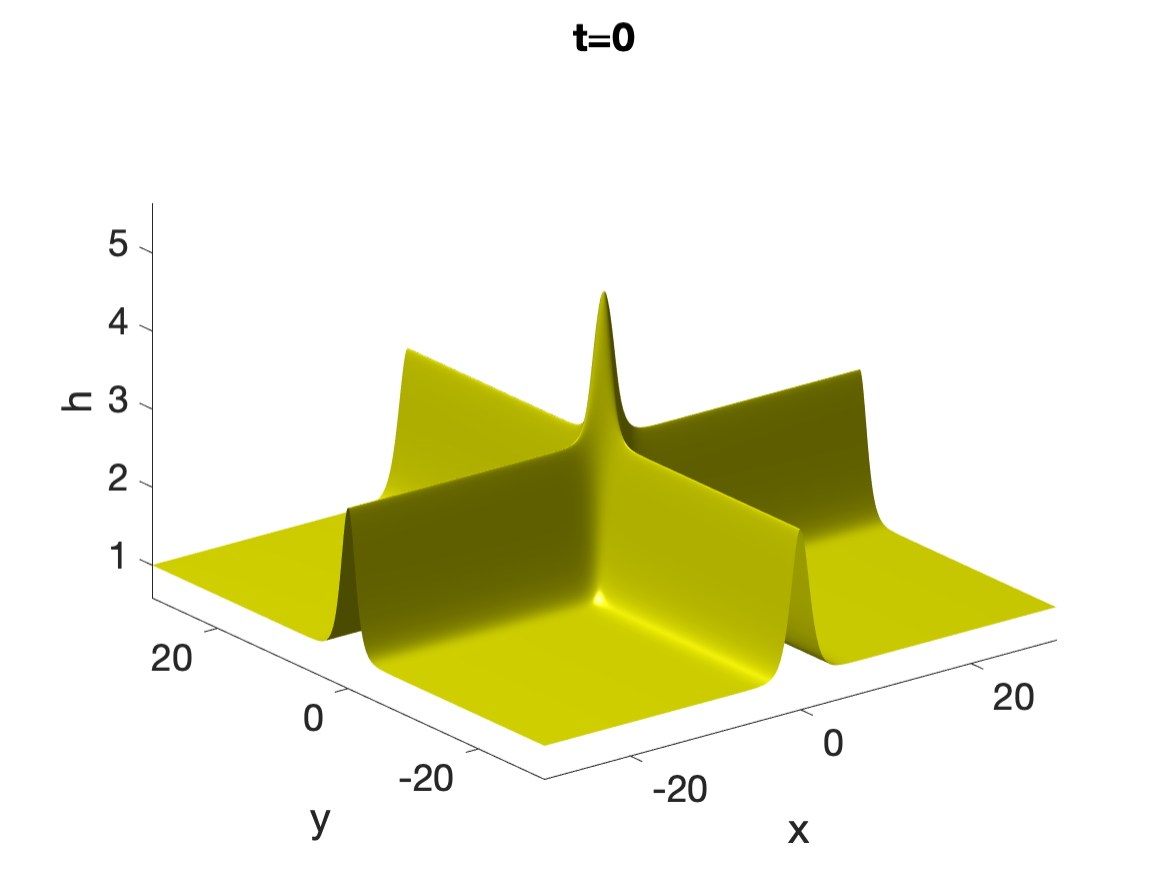}
\includegraphics[width=0.49\hsize]{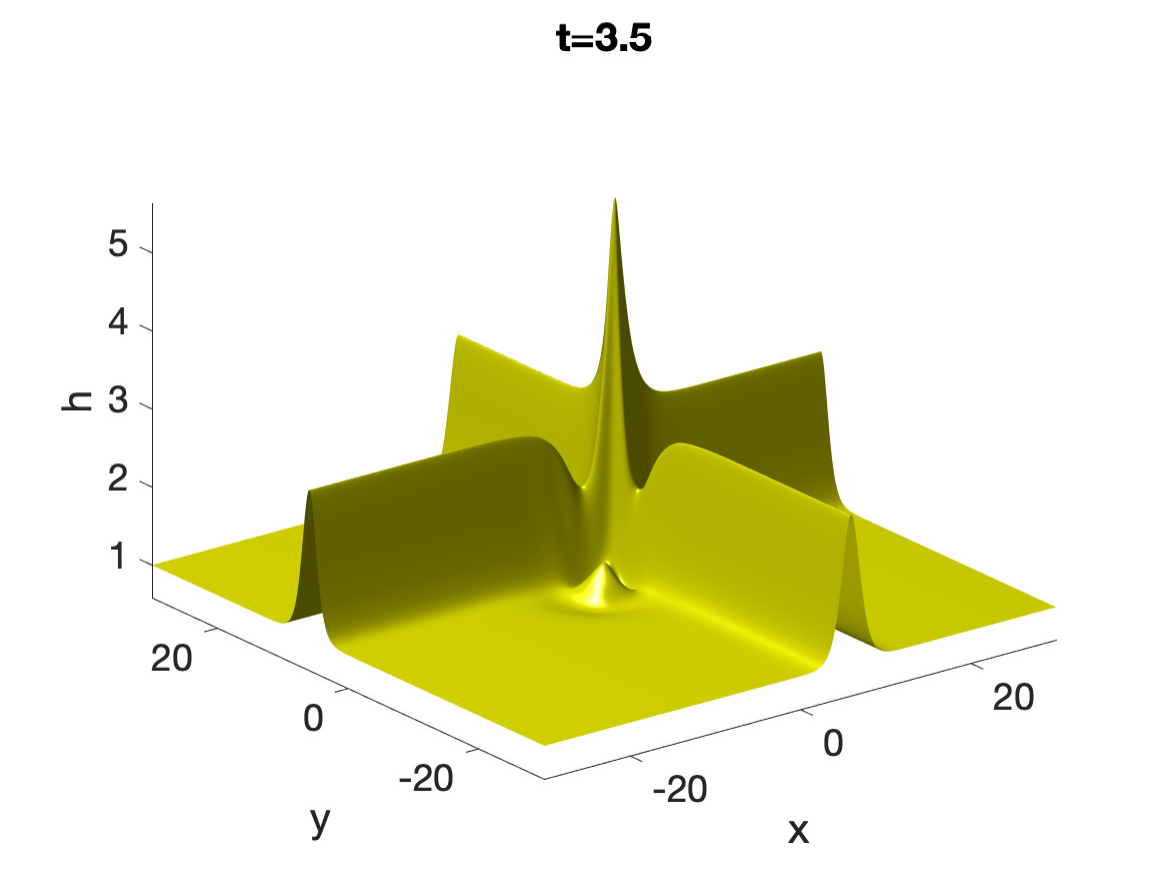}\\
\includegraphics[width=0.49\hsize]{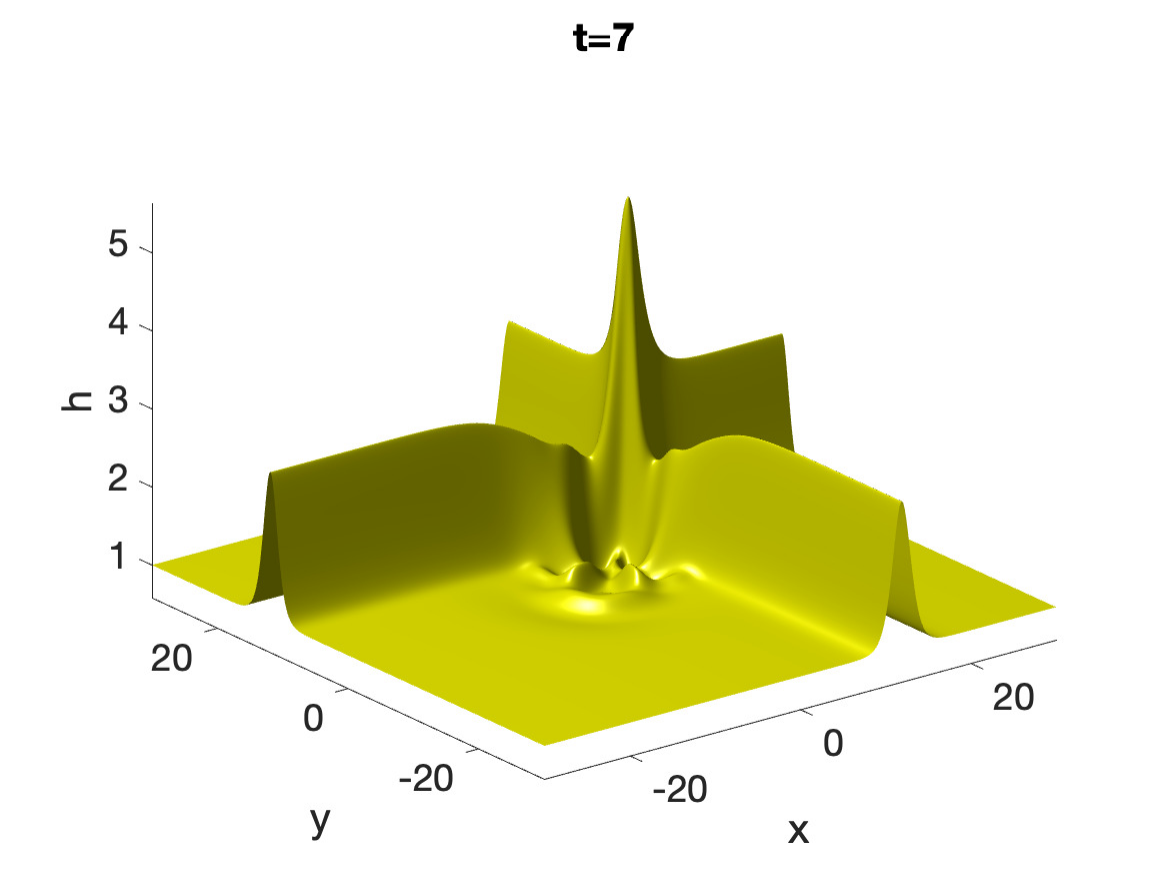}
\includegraphics[width=0.49\hsize]{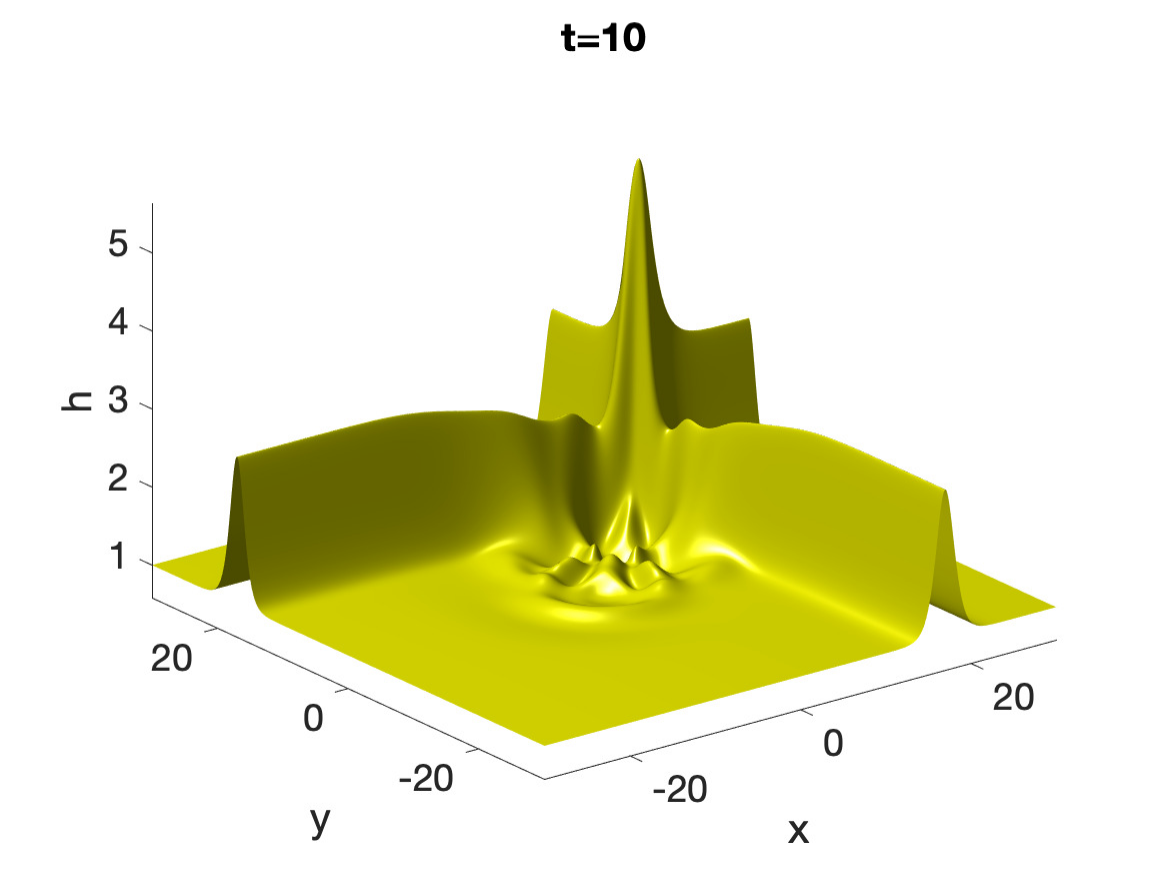}
\caption{Solution $h$ to the 2D SGN equation for initial data of the 
form (\ref{cross}) for several times. }
\label{crossh}
\end{figure}

The solution $h$ for the initial data (\ref{cross}) is shown for 
several values of time in Fig.~\ref{crossh}. Note that this is not a perturbation of one line solitary wave, but simply 
a formal sum of two such solutions. Since the SGN equation is 
nonlinear, no dynamics close to a single line solitary wave can be 
expected in such a case.
It can be seen that the solitary waves decouple from the central 
hump. Even at longer times, we do not discover a stable structure 
localised in 2D. 

\begin{figure}[!htb]
\includegraphics[width=0.49\hsize]{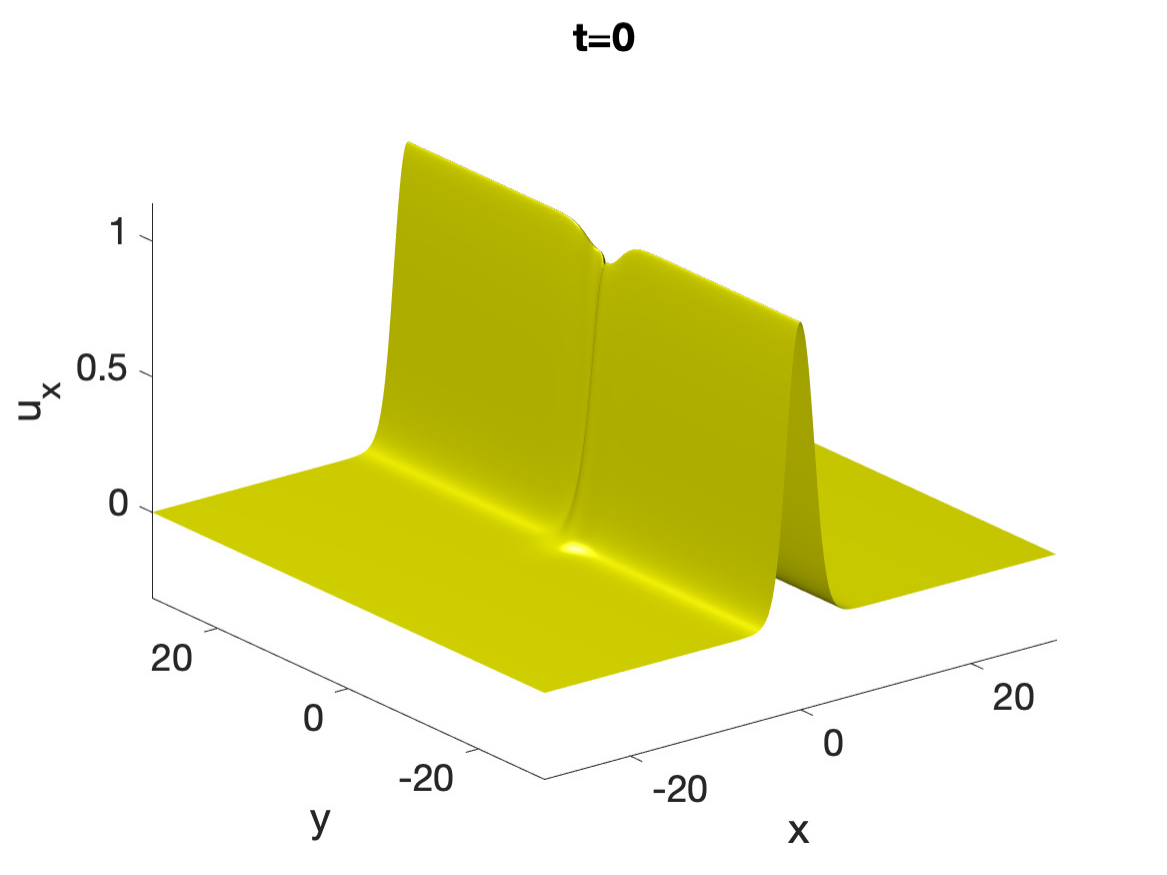}
\includegraphics[width=0.49\hsize]{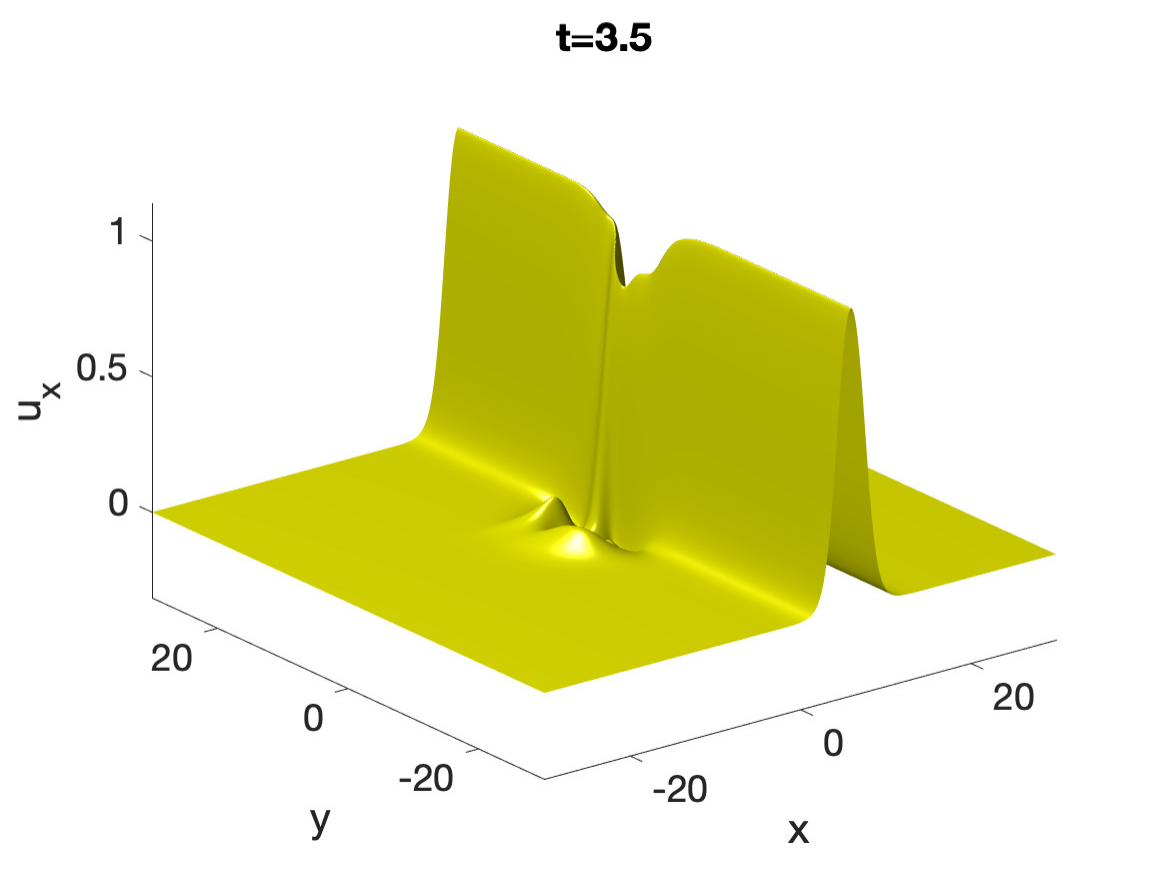}\\
\includegraphics[width=0.49\hsize]{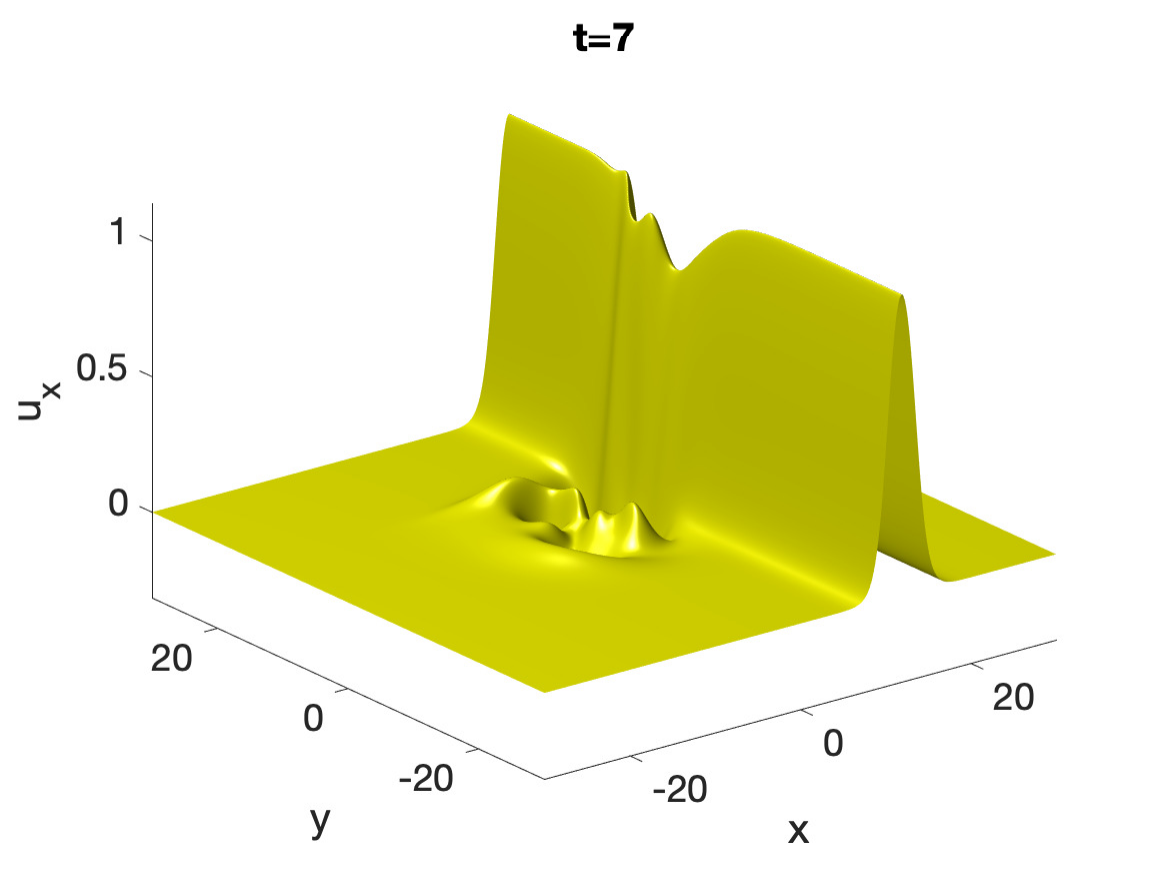}
\includegraphics[width=0.49\hsize]{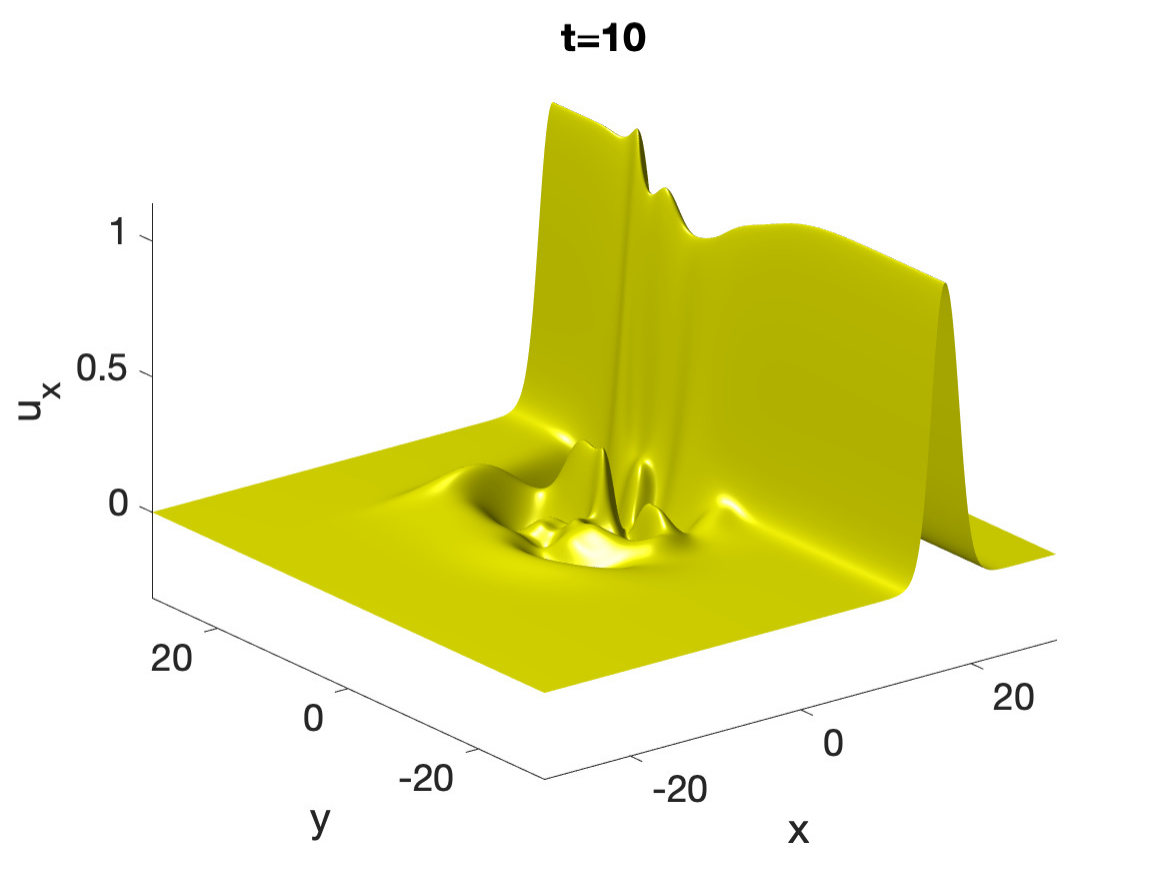}
\caption{Solution $u_{x}$ to the 2D SGN equation for initial data of the 
form (\ref{cross}) for several times. }
\label{crossux}
\end{figure}

The solution $u_{x}$ is localized in $y$-direction as can be seen in 
Fig.~\ref{crossux}. It appears to 
slowly decompose near the center where the second solitary wave is 
superposed. For symmetry reasons the $u_{y}$ solution is just the 
figure for $u_{x}$ 
rotated by 90 degrees and therefore not shown. 

\section{Localized initial data}
In this section we study  the time evolution of localized initial 
data. It is shown that initial data with radial symmetry develop into 
some annular structure with some depression near the center. It is 
argued that the solution near this center can be asymptotically 
characterized by a radially symmetric solution of the SGN 
equations. We did not find a stable structure with such initial 
data as a lump soliton for the KP I which provides support for the 
second part of the Main conjecture.

\subsection{Gaussian initial data}

As an example for localised initial data, we consider 
\begin{equation}
	h(x,y,0) = \alpha \exp(-x^{2}-y^{2}),\quad v_{x} = v_{y}=0
	\label{gauss},
\end{equation}
i.e., Gaussian initial data. Note that the code does not use the 
radial symmetry of the initial data. 

\begin{figure}[!htb]
\includegraphics[width=0.49\hsize]{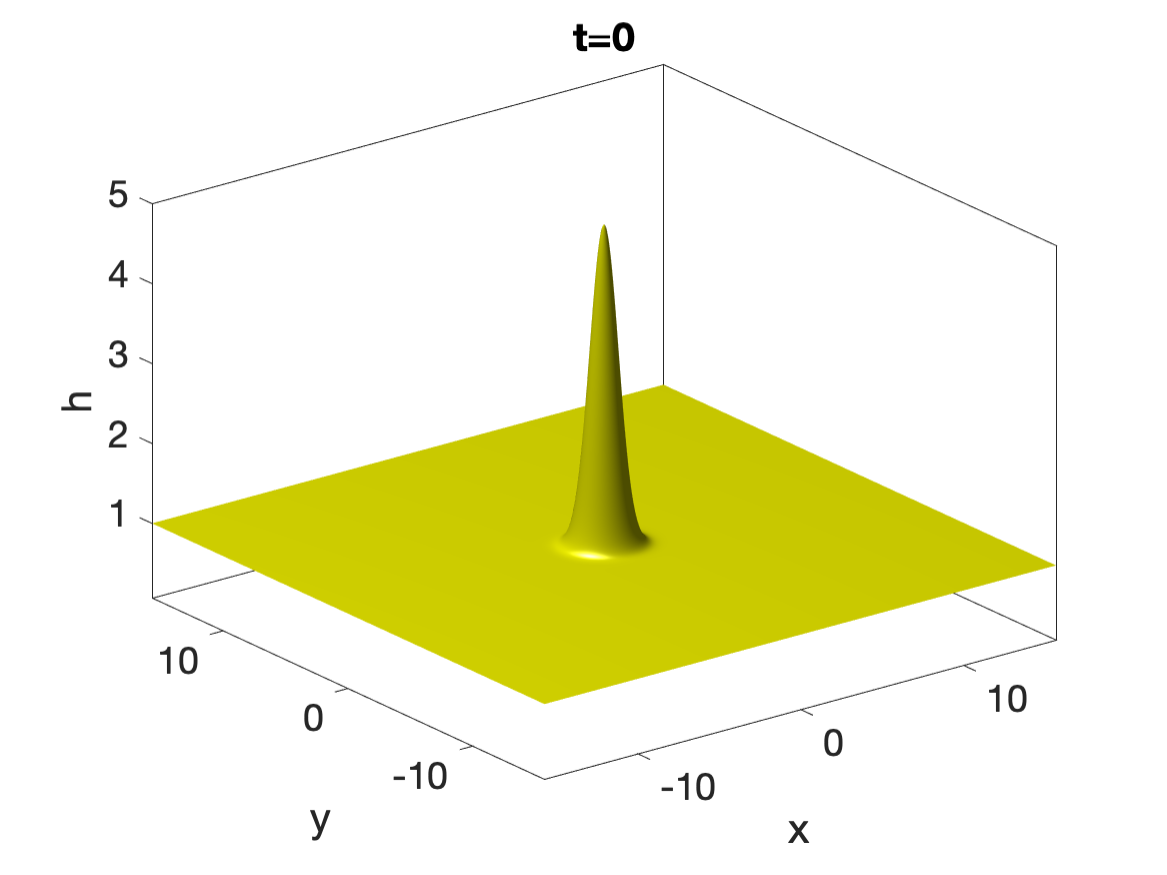}
\includegraphics[width=0.49\hsize]{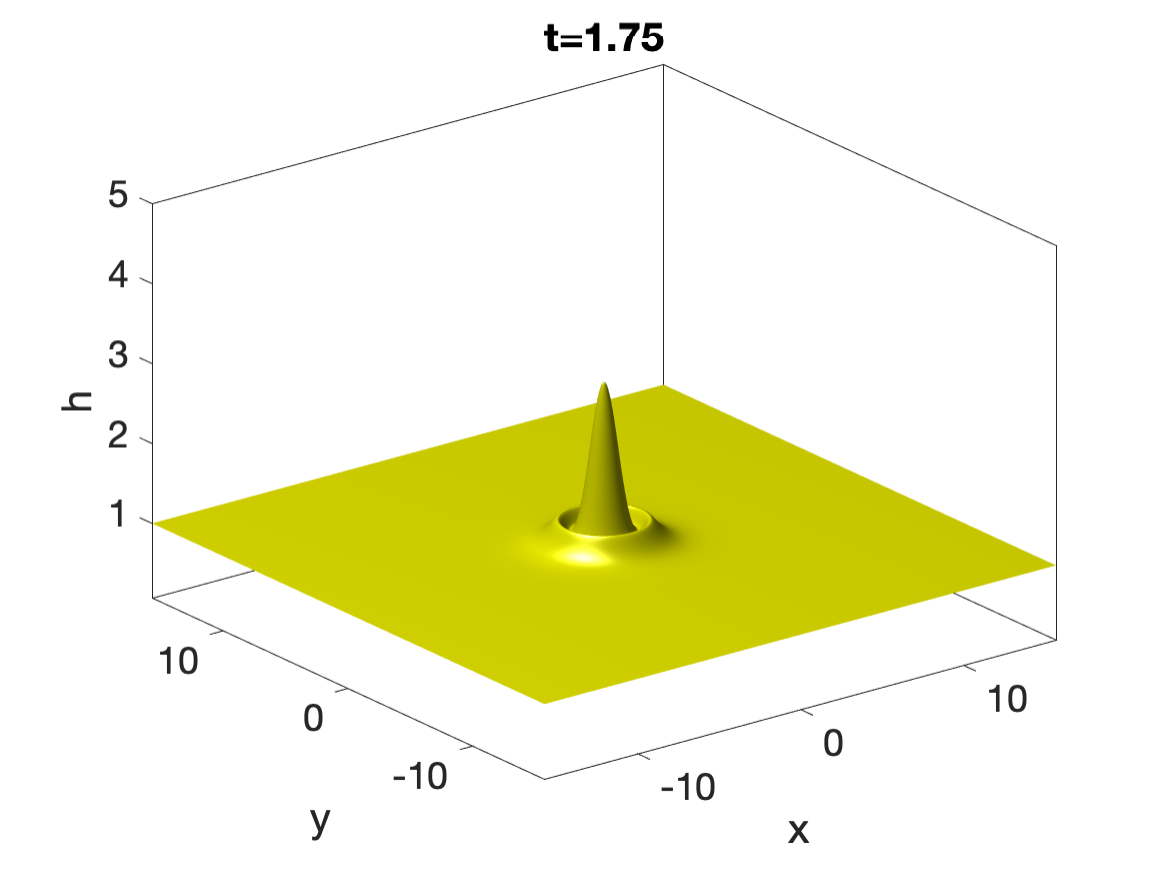}\\
\includegraphics[width=0.49\hsize]{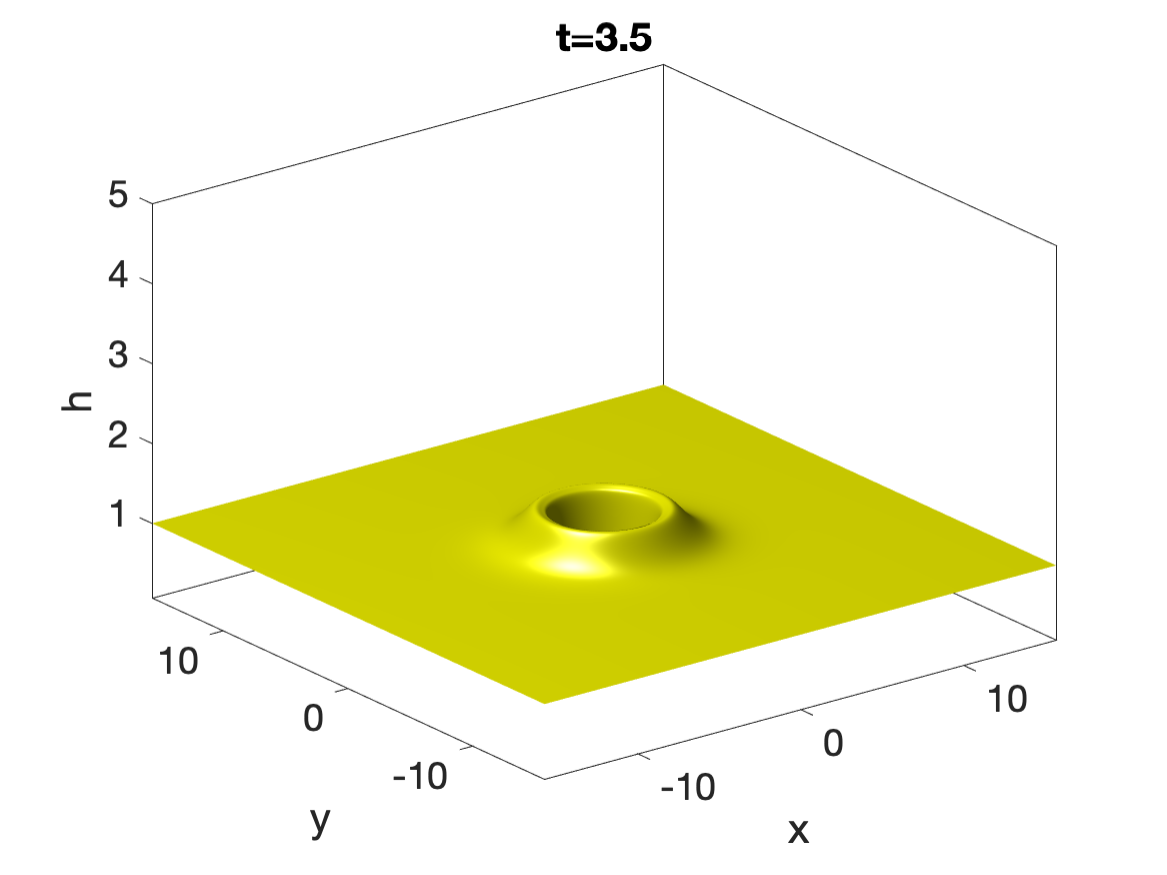}
\includegraphics[width=0.49\hsize]{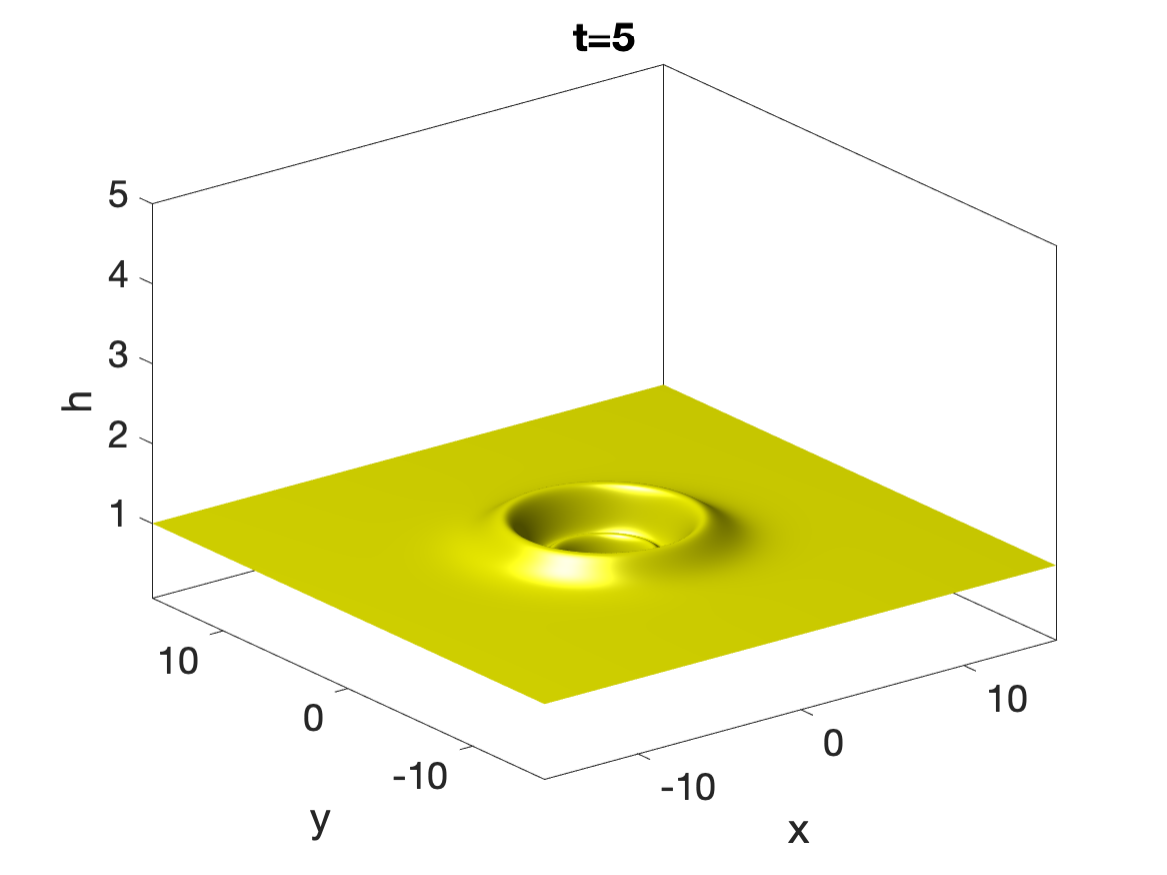}
\caption{Solution $h$ to the 2D SGN equation for initial data of the 
form (\ref{gauss}) for several times. }
\label{gaussh}
\end{figure}

For $\alpha=4$, we use 
$N_{x}=N_{y}=2^{10}$ Fourier modes and $L_{x}=L_{y}=5$ with 
$N_{t}=10^{3}$ time steps for $t\leq 5$. Relative energy conservation 
is of the order of $10^{-10}$. The solution $h$ for several 
times can be seen in Fig.~\ref{gaussh}. The initial hump evolves into 
some annular structure with radius increasing with time. Very similar 
figures can be obtained for $\alpha=2$ or $\alpha=8$. The SGN 
equation appears to be defocusing in the sense that localised humps 
get dispersed. 

\begin{figure}[!htb]
\includegraphics[width=0.49\hsize]{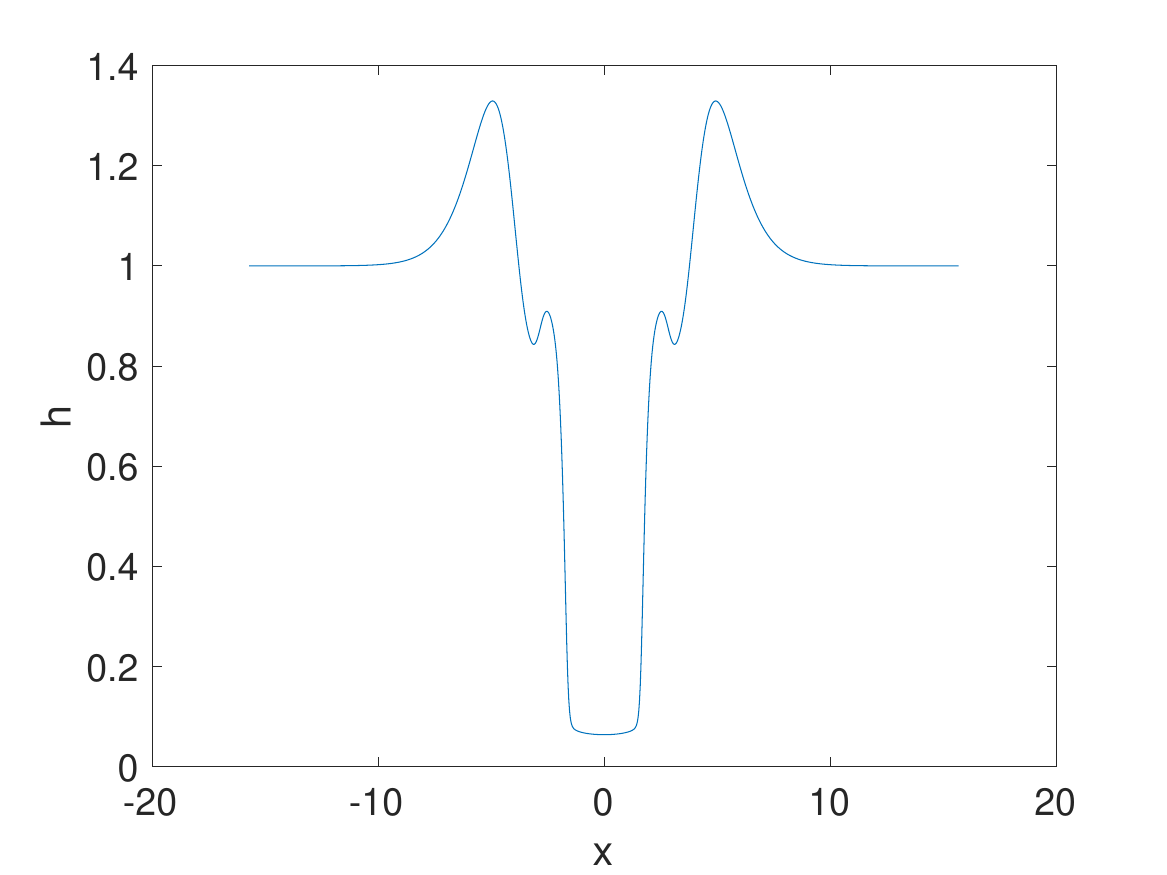}
\includegraphics[width=0.49\hsize]{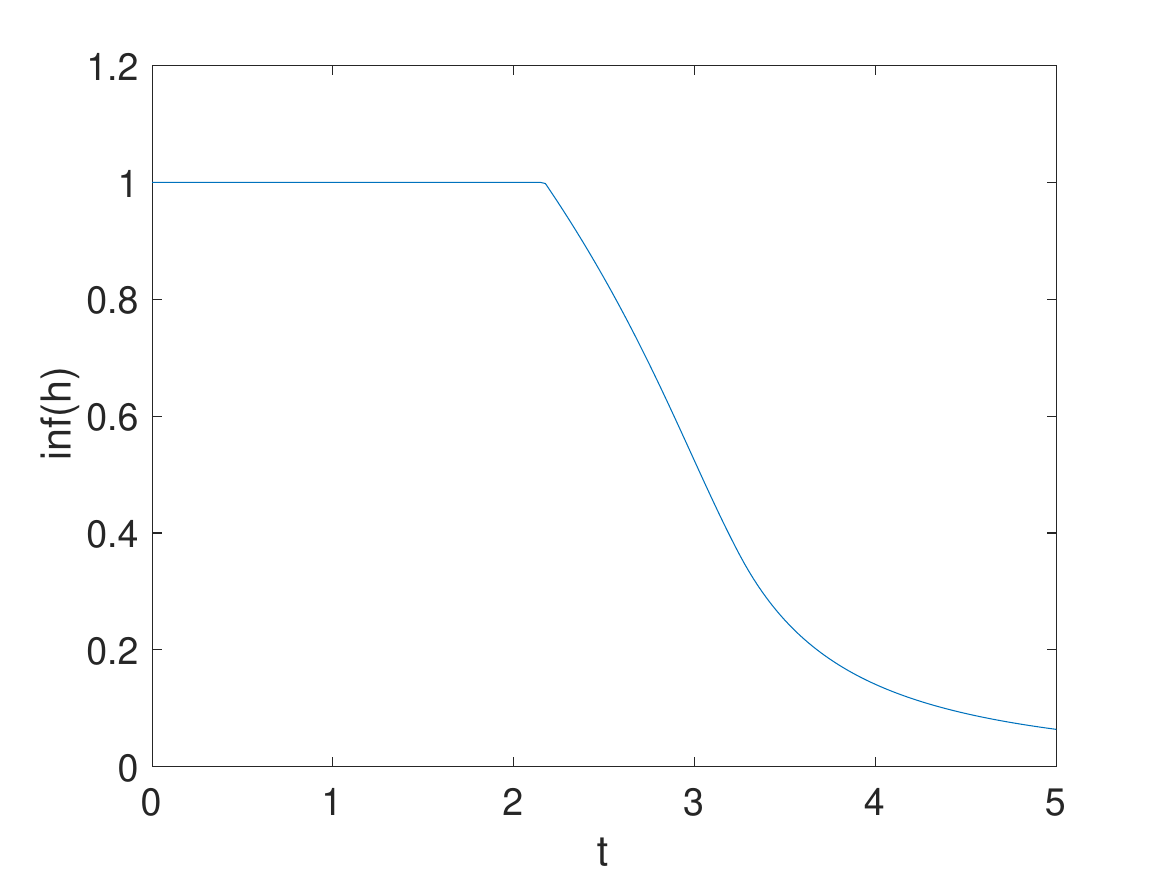}
\caption{Solution $h$ to the 2D SGN equation for initial data of the 
form (\ref{gauss}) for $t=5$ on the $x$-axis on the left, and the 
infimum of the solution in dependence of time on the right. }
\label{gaussmin}
\end{figure}

Interestingly the initial hump turns into a depression wave in finite 
time, a wave with smaller elevation than the asymptotic value for 
$h$. To illustrate this better, we show in Fig.~\ref{gaussmin} on the 
left the solution $h$ at the last shown time in Fig.~\ref{gaussh} on 
the $x$-axis. It can be seen that the solution is getting close to the 
bottom of the basin. The infimum of the solution in dependence of 
time is shown on the right of the same figure. At $t\sim 2.5$, this 
infimum becomes smaller than $h_{\infty}$ and continues to decrease. 
The no-cavitation condition appears to be always satisfied, and the 
DFT coefficients indicate that the solution stays smooth for all 
considered times. 

\begin{figure}[!htb]
\includegraphics[width=0.49\hsize]{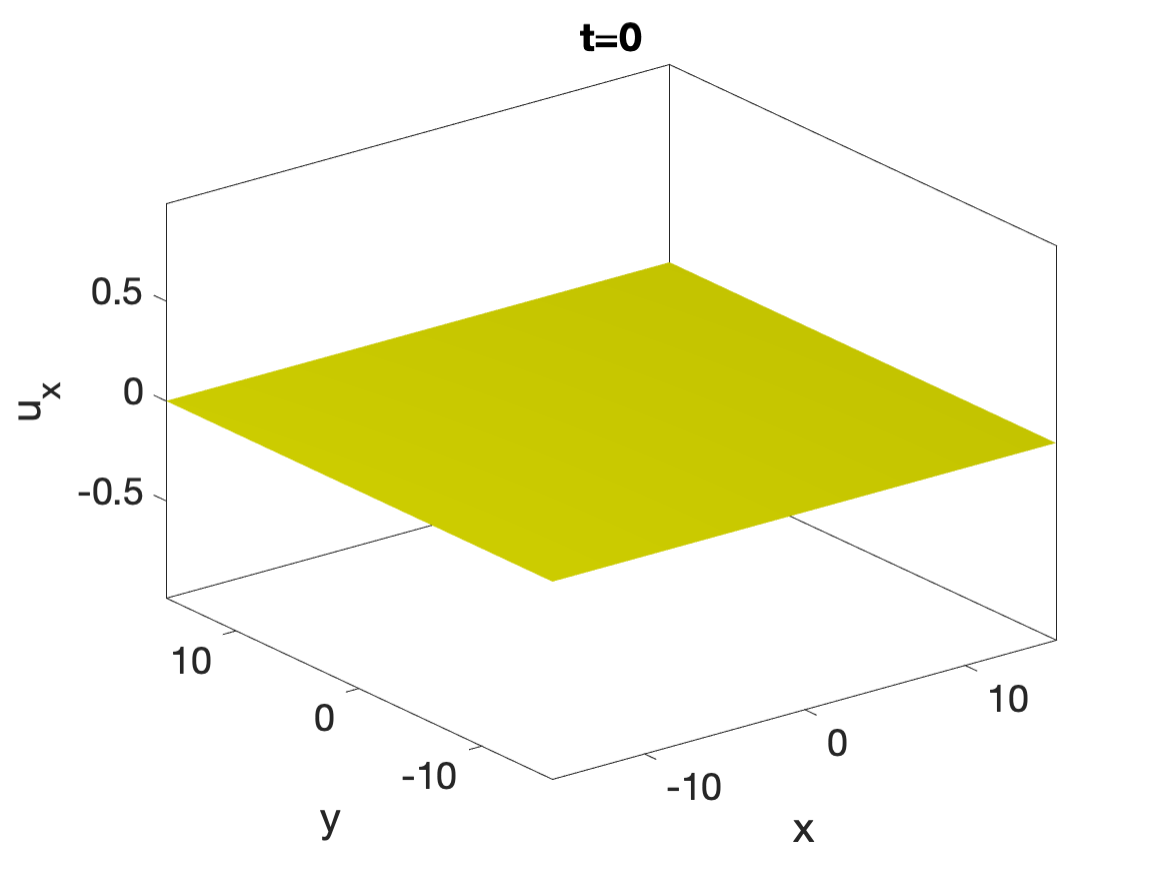}
\includegraphics[width=0.49\hsize]{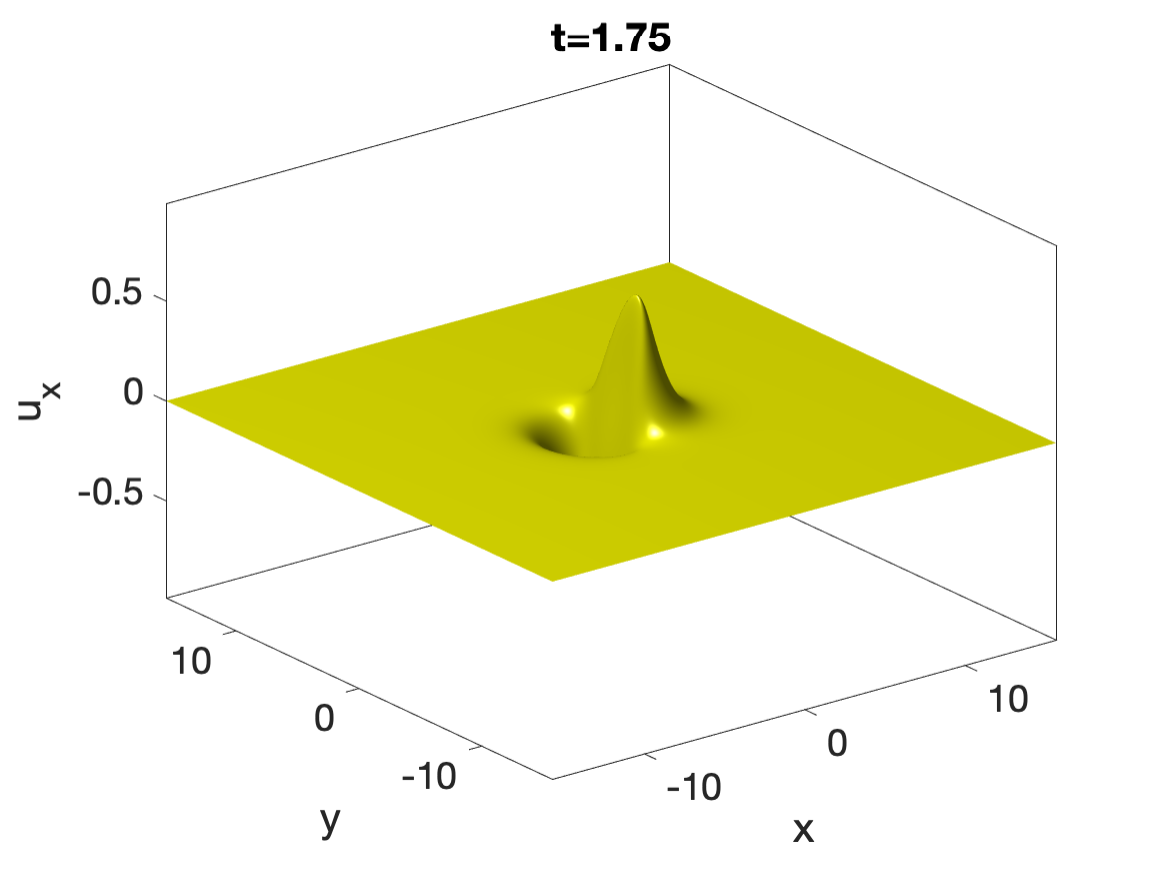}\\
\includegraphics[width=0.49\hsize]{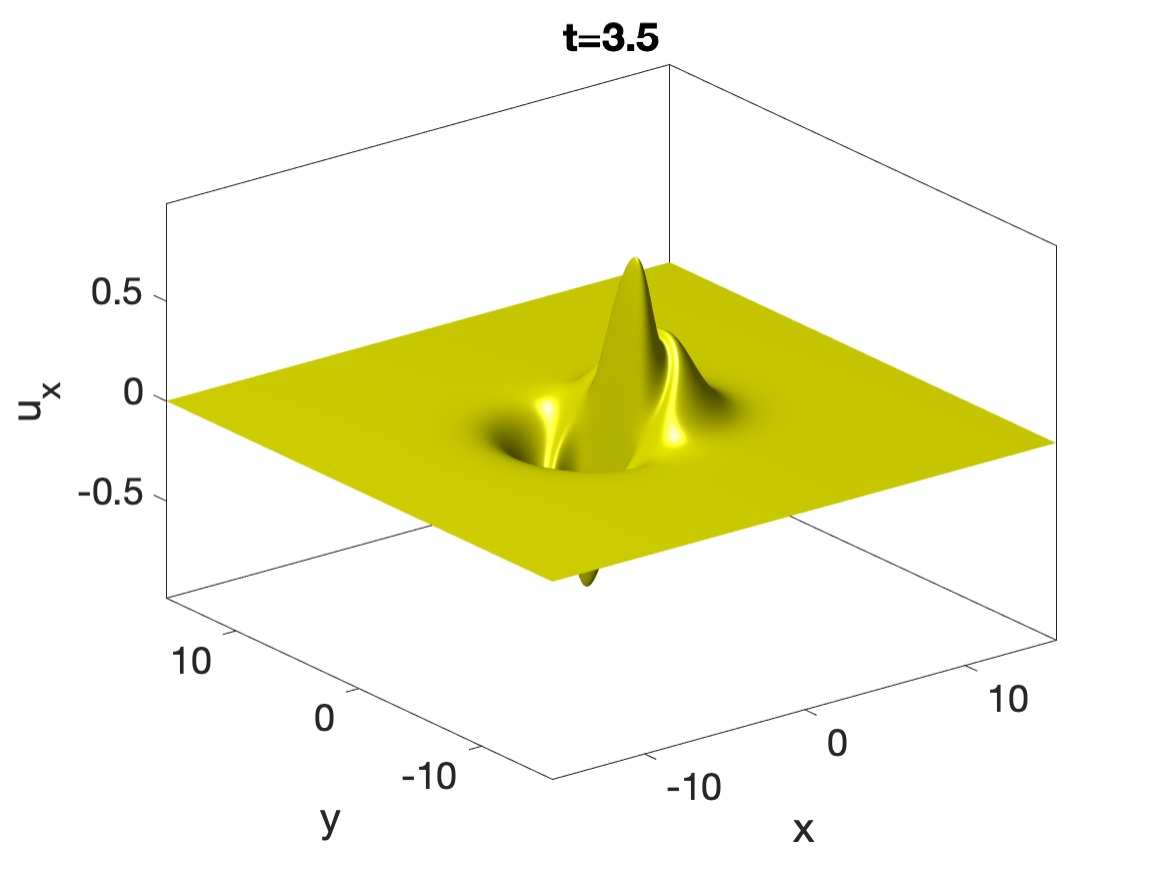}
\includegraphics[width=0.49\hsize]{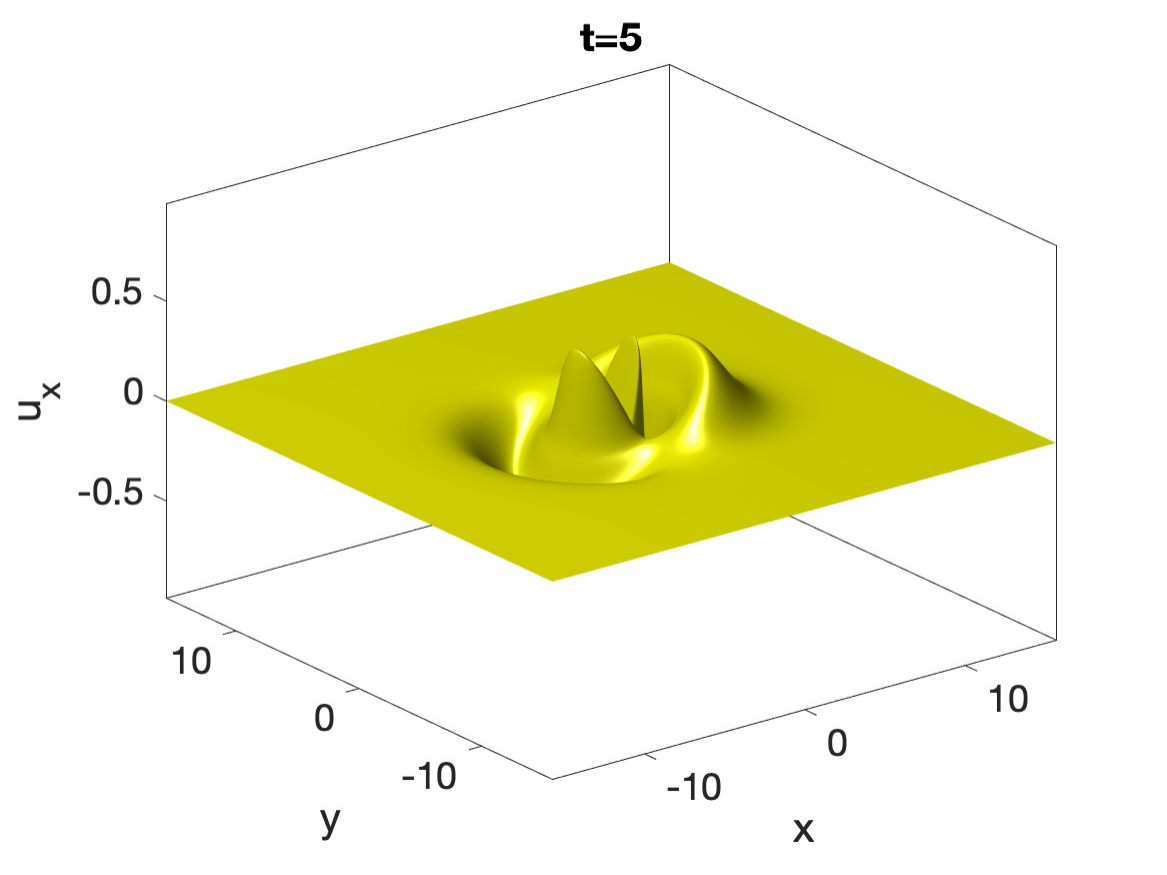}
\caption{Solution $u_{x}$ to the 2D SGN equation for initial data of the 
form (\ref{gauss}) for several times. }
\label{gaussux}
\end{figure}

The solution $u_{x}$ to the SGN equation for the initial data 
(\ref{gauss}) can be seen for several values of $t$ in 
Fig.~\ref{gaussux}. The function $u_{y}$ shows the same behavior 
rotated by 90 degrees for symmetry reasons and is not shown. 

\subsection{Radially symmetric solutions}

\begin{figure}[!htb]
\includegraphics[width=0.49\hsize]{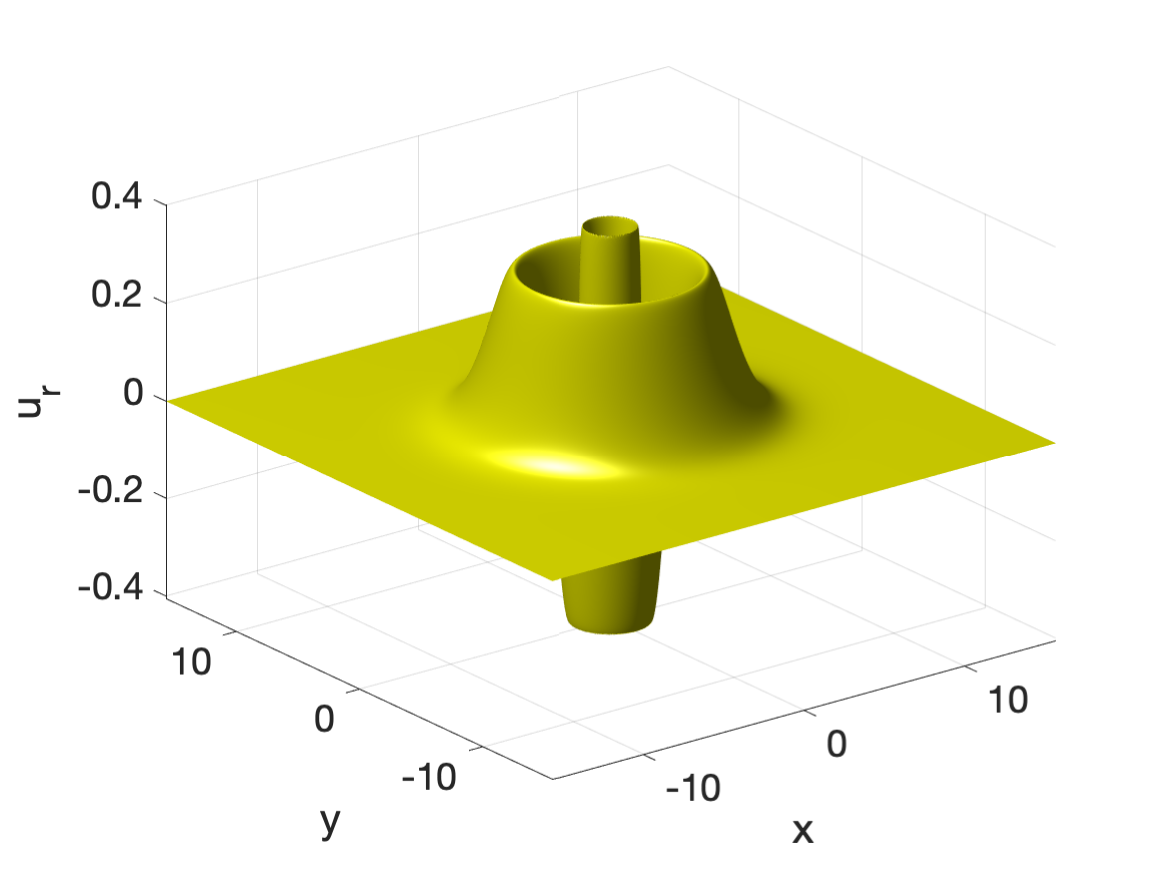}
\includegraphics[width=0.49\hsize]{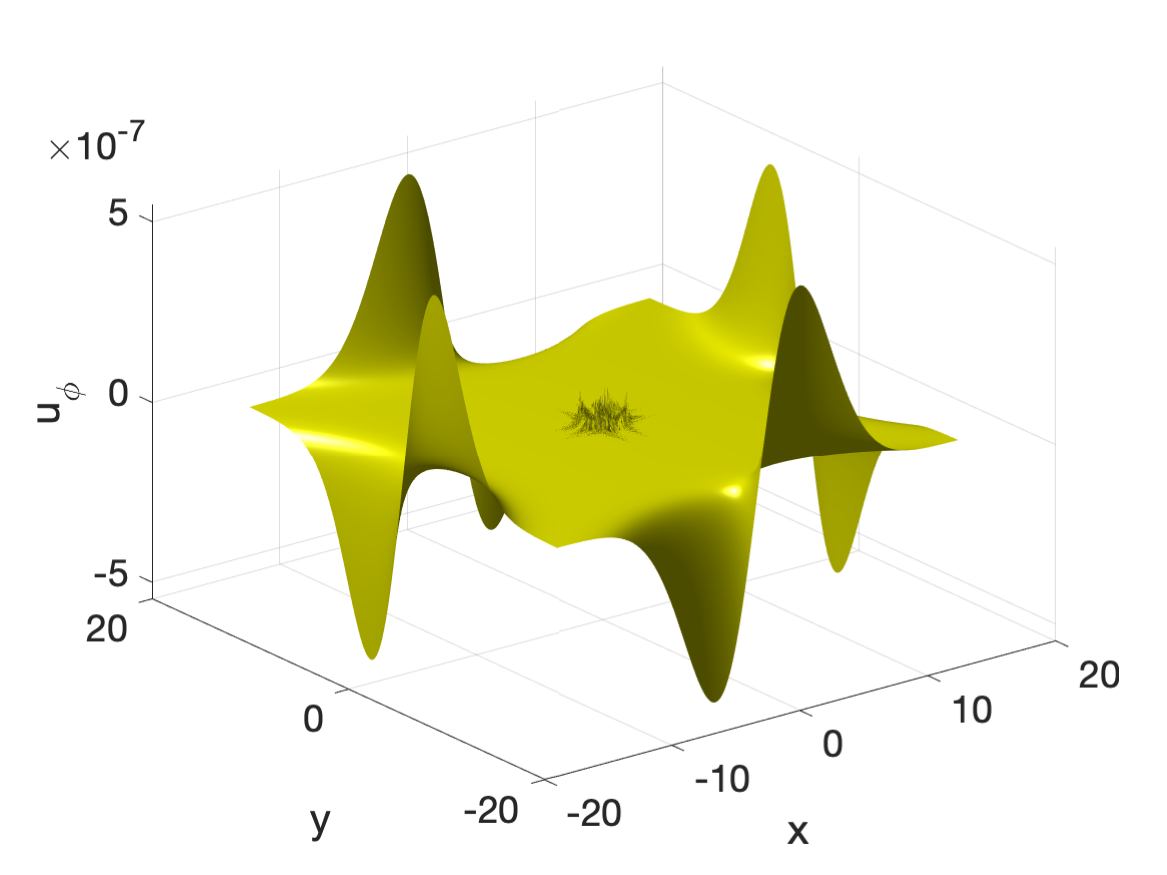}
\caption{Velocity $\overline{\mathbf{u}}$ for the 2D SGN equation for initial data of the 
form (\ref{gauss}) for $t=5$, on the left the radial component, on 
the right the azimuthal component. }
\label{gaussupolar}
\end{figure}

	The formation of a ``cylinder-type" structure  in the case of 
	axisymmetric localized initial data  suggests that the solution 
	for long times becomes radially symmetric	near the center of 
	the annular structure: the position of the free 
	surface  is almost flat, no obvious oscillations can be seen. 
In the radially symmetric case, we can construct  an exact solution 
to  the SGN equations approximating this ``flat'' structure. More 
exactly, consider the SGN equations \eqref{eq:SGN} in the  polar coordinates 
($x=r\cos \phi$, $y=r\sin\phi$). In the case where the azimuthal velocity $u_\phi$  is zero and all other variables depend only on $r$, the governing equations  for  $h$ and radial velocity $u_r$ take a standard form:
	\begin{equation*}
	\frac{\partial (rh)}{\partial t} +\frac{\partial (rhu_r)}{\partial r}=0, \quad h\left(	\frac{\partial u_r}{\partial t}+u_r	\frac{\partial u_r}{\partial r}\right)+	\frac{\partial p}{\partial r}=0, \quad r=\sqrt{x^2+y^2}
   	\end{equation*}
They admit the following exact solution 
\begin{equation}
	\label{profile}
\tilde h(t)=\frac{h_0}{(c_0+w_0 t)^2}, \quad \tilde u_r(t,r)=\frac{w_0 r}{c_0+w_0 t},
\end{equation}
with constants $c_0$, $h_0$ and $w_0$.
We show the radial and azimuthal component of the velocity for $t=5$ 
for the situation shown in Fig.~\ref{gaussux} in 
Fig.~\ref{gaussupolar}. It can be seen that the radial velocity shows 
almost the expected radial symmetry whereas the azimuthal component 
is of the order of $10^{-7}$ (still much larger than the estimated 
numerical error) and almost vanishing near the origin. 

\begin{figure}[!htb]
\includegraphics[width=0.49\hsize]{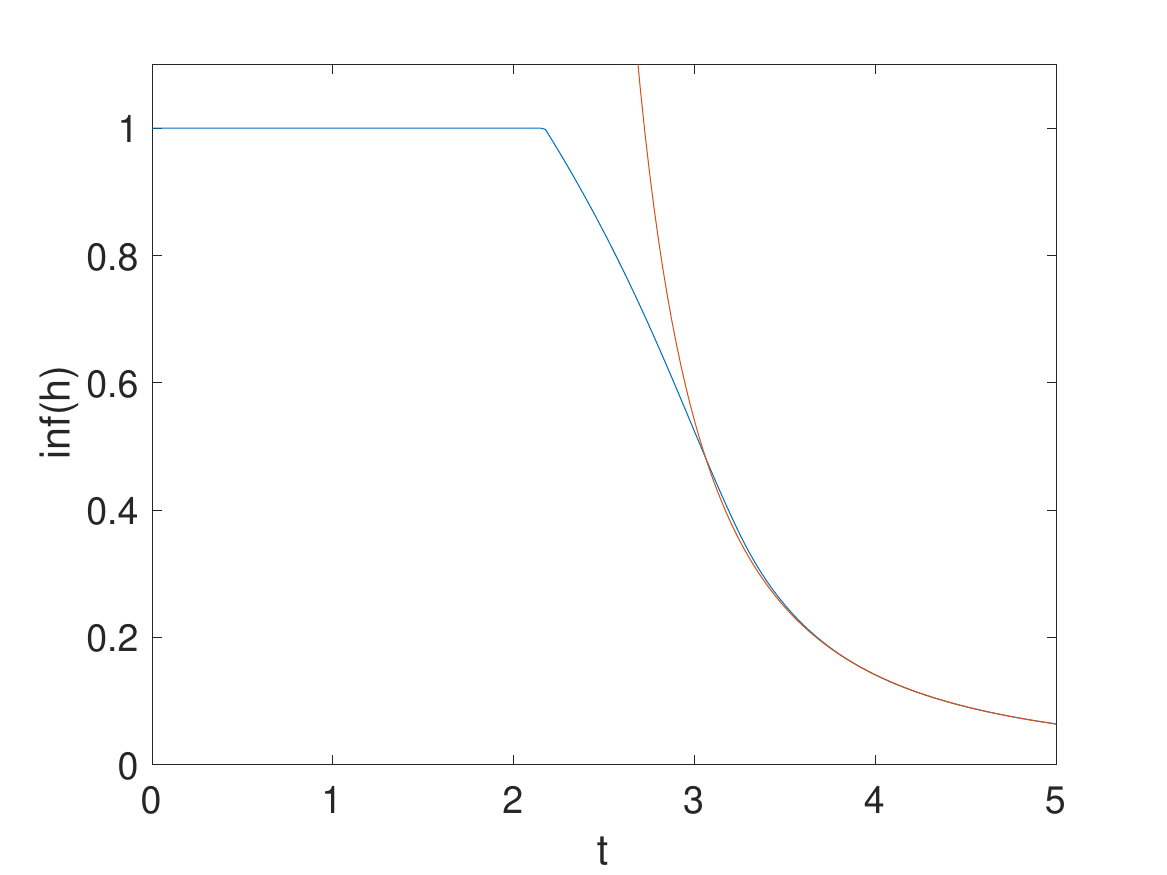}
\includegraphics[width=0.49\hsize]{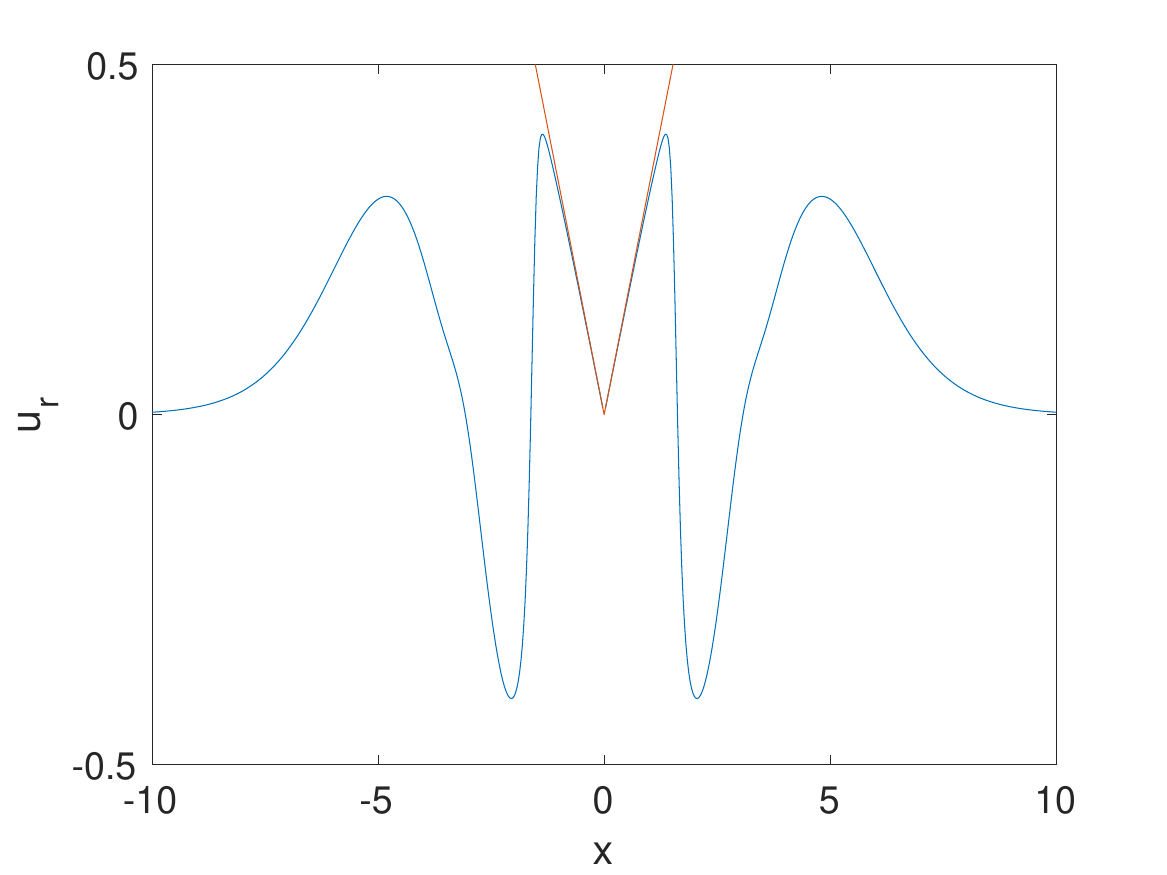}
\caption{Infimum of the solution $h$ to the 2D SGN equation for initial data of the 
form (\ref{gauss}) with the fitted profile (\ref{profile}) in red  on the left, and the 
and the radial component of the velocity in dependence of $x$ with 
the fitted solution (\ref{profile}) in red  on the right. }
\label{gaussminfit}
\end{figure}

We can fit the infimum of the solution $h$ in Fig.~\ref{gaussmin} on 
the left to the asymptotic formula. This means we choose the 
parameters $h_{0}/c_{0}^{2}$ and $w_{0}/c_{0}$ such that 
$h(1+tw_{0}/c_{0})^{2}$ is roughly constant. To this end we apply the 
algorithm \cite{fminsearch} being distributed with Matlab as the 
command \textit{fminsearch} for $t\geq 3.75$ (the result does not 
change much if this time is slightly varied). We get 
$h_{0}/c_{0}^{2}=0.1545$ and $w_{0}/c_{0}=-0.5117$. The result of the 
fitting of the infimum of the solution can be seen on the left of 
Fig.~\ref{gaussminfit}. It can be seen that the fitting  gives an 
excellent result (the residual on the fitted domain is of the 
order of $10^{-3}$) which gives a strong indication that the solution 
near the origin follows the dynamics of the special solution \eqref{profile}. This is 
further confirmed by comparing the solution  $ \tilde u_r$ with the fitted 
parameters to the radial component $u_r$ of the velocity as shown on 
the right. The linear behavior in $r$  can be clearly 
recognized. 

\section{Conclusion}
In this paper we have presented a numerical approach to the 2D SGN 
equations based on a Fourier spectral method with the Krylov subspace 
technique GMRES to solve an elliptic equation. The code allowed to 
provide strong numerical evidence for the transverse stability of the line solitary waves of SGN 
as well as the defocusing character of the SGN equations in 2D. No 
stable structures localised in 2D were discovered. It was shown that 
radially symmetric initial data evolve into an annular structure with 
a depression near the center that can be asymptotically described via 
a special radially symmetric solution to the SGN equations. The 
no-cavitation condition appears to be satisfied at all times. 

Stronger computers than the ones we had access  for this paper 
(essentially a laptop) will 
be needed to address questions having been studied in 1D in 
\cite{DK}. An interesting such question will be the appearence of 
\textit{dispersive shock waves} (DSW), zones of rapid modulated 
oscillations in the solutions. They appear in dispersive PDEs in the 
vicinity of shocks of the solutions to the corresponding 
dispersionless PDEs for the same initial data \cite{El_Hoefer_2016}. 
While 1D DSW for the SGN equations have already been actively 
investigated \cite{El06, Gavrilyuk_2020, Pitt_2018}, their 
multi--dimensional analogues are still waiting to be studied. 


Another interesting question is whether the SGN equation allows for 
shocks in finite time for smooth initial data. A comparison to the 
Camassa-Holm (CH) equation was presented in 1D in \cite{DK}.  Since 
both CH and SGN are nonlocal and since it is 
known that CH solutions can have a blow-up (an explosion of some norm 
of the solution) for certain classes of initial data,  a similar 
behaviour is possible for the SGN equations. 
%

This will be the subject of 
further research.


\begin{thebibliography}{99}
%
%
%
%
%
%
%
%
%
%
%
%
	
	
	\bibitem{DK} V. Duch\^ene, C. Klein, Numerical study of the 
	 Serre-Green-Naghdi equations and a fully dispersive 
	 counterpart, 	Discr. \& Cont. Dyn. Syst.  B, 27(10), 5905-593 (2021) doi: 
	 10.3934/dcdsb.2021300   
	 
	 \bibitem{El06} G. A.  El,  R. H. J. Grimshaw and N. F. Smyth  2006 
	 Unsteady undular bores in fully nonlinear shallow-water theory, Phys.
	 Fluids \textbf{18}  027104.
	 
	 \bibitem{El_Hoefer_2016} G. A. El and M. A. Hoefer, Dispersive shock waves and modulation theory, Physica D 333, 11-65 (2016)
	 
	\bibitem{Gavrilyuk_2014} S.  Gavrilyuk, H.  Kalisch \&  Z.  Khorsand, A kinematic conservation law in free surface flow, Nonlinearity, {\bf 28} (2014) 1805--1821.
	
	
	\bibitem{Gavrilyuk_Teshukov_2001} S.  L.  Gavrilyuk \&  V. M.  Teshukov, 	Generalized vorticity for bubbly liquid and dispersive shallow water equations, Continuum Mechanics and Thermodynamics, {\bf 13} (2001) 365-382
	

%
%
	\bibitem{Gavrilyuk_2020} S. Gavrilyuk, B. Nkonga,  K.--M. Shyue \&
	L. Truskinovsky,
	Stationary  shock-like transition fronts 
	in dispersive systems, Nonlinearity, \textbf{33} (2020),
	5477-5509.
	
	
	\bibitem{Green_74} A. E. Green, N. Laws \& P. M. Naghdi, On the
	theory of water waves, Proc. R. Soc. Lond. \textbf{A 338} (1974), 43--55. 
	
	\bibitem{Green_76} A. E. Green \& P. M. Naghdi, A derivation
	of equations for wave propagation in water of variable depth, J. Fluid
	Mech. \textbf{78} (1976), 237--246.
	
\bibitem{etna} C.~Klein,  Fourth order time-stepping for low dispersion Korteweg-de 
Vries and nonlinear Schr\"odinger equation,  ETNA Vol. 29 116-135 (2008).

	\bibitem{KSbook} C.~Klein and J.-C.~Saut, Nonlinear dispersive equations --- 
	Inverse Scattering and PDE methods, 
	Applied Mathematical Sciences 209 (Springer, 2022)
	
	 \bibitem{Kod}Y. Kodama, KP Solitons and the Grassmannians: 
	 Combinatorics and Geometry of Two-Dimensional Wave Patterns, 
	 SpringerBriefs in Mathematical Physics \textbf{22} (2017).

	 \bibitem{Kra} R. Krasny, A study of singularity formation in a 
	 vortex sheet by the point-vortex
approximation, J. Fluid Mech. 167 (1986) 65–93.
	
	\bibitem{LannesBOOK_2013}  D. Lannes, The Water Waves Problem, 
	Mathematical Surveys and Monographs, vol. \textbf{188} (Amer. Math. Soc.,
	Providence, 2013). 
	
	
	\bibitem{Makarenko_1986}  N. Makarenko, A second long-wave approximation
	in the Cauchy-Poisson problem, Dynamics of Continuous Media, v. 77 (1986), pp. 56-72
	(in Russian). 
	
\bibitem{fminsearch}  J. C. Lagarias, J. A. Reeds, M. H. Wright, and 
P. E. Wright. Convergence Properties of the Nelder-Mead Simplex 
Method in Low Dimensions. SIAM Journal of Optimization. Vol. 9, 
Number 1, 1998, pp. 112-147. 

	\bibitem{Metayer10} O. Le M\'etayer, S. Gavrilyuk \& S. Hank, 
	A numerical scheme for the Green-Naghdi model, J. Comp. Phys. \textbf{229} (2010),
	2034--2045. 
	
	
\bibitem{Li_2001} Yi A. Li, Linear stability of solitary waves of the Green-Naghdi	Equations, Communications on Pure and Appl. Math.  \textbf{LIV} (2001), 501--536. 
	
%
%
	
	
%
%
	
%
	
		\bibitem{Pitt_2018} J. P. A. Pitt, C. Zoppou \& S. G. Roberts, 
		Behaviour of the {S}erre equations in the presence of
		steep gradients revisited, Wave Motion 76 (2018), 61--77
		
	
	\bibitem{Salmon} J. Miles \& R. Salmon, Weakly dispersive
	nonlinear gravity waves, J. Fluid Mechanics \textbf{157} (1985), 519--531. 
	
	\bibitem{GMRES} Y. Saad and M. H. Schultz. GMRES: a generalized minimal residual algorithm for solving nonsymmetric linear systems. SIAM J. Sci. Statist. Comput., 7(3):856–869, 1986.


	
	\bibitem{Salmon_1998} R. Salmon,   Lectures on Geophysical Fluid Dynamics,  Oxford University Press, 1998,  ISBN 9780195355321.
	
	
	\bibitem{Serre_53} F. Serre, Contribution \`a l'\'etude des \'ecoulements permanents
	et variables dans les canaux, La Houille Blanche	\textbf{8} (1953), 374--388. 
	


	
	\bibitem{Su_Gardner_1969} C. H. Su \& C. S.  Gardner,  Korteweg - de Vries Equation
	and Generalisations. III. Derivation of the Korteweg - de Vries Equation and Burgers
	Equation,  J.  Math. Physics,  \textbf{10} (1969) 536--539. 
	
	
	\bibitem{Tkachenko_2023}	S. Tkachenko, S. Gavrilyuk \& J. Massoni, Extended Lagrangian approach for the numerical study of multidimensional dispersive waves: applications to the Serre-Green-Naghdi equations, Journal of Computational Physics (2023), 111901.
	
	
	
	
	\bibitem{trefethen} L. Trefethen. Spectral Methods in Matlab. SIAM, Philadelphia, PA, 2000.
\end{thebibliography}
\end{document}